\documentclass[12 pt]{report}

\usepackage{amsthm}
\usepackage{amssymb}
\usepackage{latexsym}
\usepackage{url}
\usepackage{graphicx}
\usepackage{epsf}
\usepackage{amsmath}
\usepackage[normal]{caption2}

\author{Orr Moshe Shalit}

\title{Guided Dynamical Systems and Applications to Functional and Partial Differential Equations}

\newtheorem{thm}{Theorem}[section]
\newtheorem{cor}[thm]{Corollary}
\newtheorem{lem}[thm]{Lemma}
\newtheorem{prop}[thm]{Proposition}

\newtheorem{conj}[thm]{Conjecture}
\newtheorem{defin}[thm]{Definition}
\newtheorem{defins}[thm]{Definitions}
\newtheorem{metadefin}[thm]{Meta-Definition}
\newtheorem{expl}[thm]{Example}

\newtheorem{rem}[thm]{Remark}

\begin{document}
\addtocontents{toc}{\protect\thispagestyle{empty}}

\thispagestyle{empty} \vspace{2cm} {\centering \Huge Guided
Dynamical Systems and Applications to Functional and Partial
Differential Equations
\par}\vspace{8cm} {\centering \huge Orr Moshe Shalit
\par}

\begin{titlepage}
\vspace{2cm} \centering
\Large Guided Dynamical Systems and Applications to Functional and Partial Differential Equations\\

\vspace{4cm} \centering \large Research Thesis \\
\large Submitted in Partial Fulfillment of the Requirements for
the Degree of Master of Science in Mathematics

\vspace{4cm}
\large Orr Moshe Shalit\\

\vspace{3cm} \large Submitted to the Senate of the Technion -
Israel
Institute of Technology \\

\vspace{1cm} \large Sivan, 5765 \quad Haifa \quad June 2005
\end{titlepage}


\thispagestyle{empty} \vspace{1cm} {\centering \large This
Research Thesis Was Done Under The Supervision of Professor Boris
Paneah in the Department of Mathematics \par} \vspace{2cm}

{\centering \large I would like to thank Prof. Paneah for
inspiring me\\
and for all that he taught me\par} \vspace{2cm}

{\centering \large I would like to thank my friend Daniel Reem for
proofreading the manuscript and for making the figures
\par} \vspace{2cm}

{\centering \large The financial support granted by the Technion
during my studies is greatly acknowledged\par} \vspace{2cm}

{\centering \large This work is dedicated to the people who
supported\\ me during my studies\\
Braha and Ilan Fabian, my in laws\\
Malka and Meir Shalit, my parents\par} \vspace{2cm}



\pagebreak

\tableofcontents \thispagestyle{empty} \pagenumbering{arabic}

\listoffigures \thispagestyle{empty} \pagebreak

\setcounter{page}{1} \thispagestyle{plain}

\vspace{2.5cm} {\centering \Large Abstract \par} \vspace{1.5cm} In
this thesis I present the concept of a \emph{guided dynamical
system}, and then I exploit this idea to solve various problems in
functional equations and partial differential equations. The
results presented here are in a sense a sequel to a series of
papers by B. Paneah published in the years 1997-2004.

The last chapter in this work is an introduction containing an
overview of this work and a comparison between known and new
results.

In the first chapter I shall first explain what a guided dynamical
system is, introducing all notations and definitions to be used in
the later chapters.

In the second chapter I will use guided dynamical systems to study
functional equations which have the form
\begin{displaymath}
f(x) - \sum_{i=1}^N a_i (x) f(\delta_i (x)) = h(x) \quad, \,\, x
\in X
\end{displaymath}
where the functions $a_i$, $\delta_i$ and $h$ are given, and $f$
is an unknown, continuous real-valued or vector-valued function
defined on (typically) a compact space $X$. For this type of
equations I present here original results regarding uniqueness and
solvability, the methods used are extensions of those introduced
by Paneah.

In the third chapter we make a detor from our main route to treat
the more esoteric problem of \emph{over-determinedness}, for which
I also present some new methods and results.

In chapter 4 I will use the results of chapters 2 and 3 to give a
necessary and sufficient condition for the unique-solvability of
the second partly characteristic boundary  value problem:
\begin{displaymath}
(m\partial_x + n\partial_y)\partial_x \partial_y u = 0 \quad {\rm
in} \quad D
\end{displaymath}
\begin{displaymath}
 u  =  g \quad {\rm on} \quad \partial D \,\,.
\end{displaymath}
In this chapter I will use Paneah's reduction of the above problem
to a Cauchy type functional equation to give the necessary and
sufficient condition in terms of the dynamical properties of a
guided dynamical system in the boundary of the problem. For
specific families of domains, this necessary and sufficient
condition is then translated to explicit conditions for the
well-posedness of this hyperbolic boundary problem.

\cleardoublepage

\pagebreak
\thispagestyle{plain} \vspace{2.5cm} {\centering \Large List of
Notations \par} \vspace{1.5cm}

\begin{tabular}{ll}

$(X,\delta)$ & dynamical system generated in $X$ by the
maps $\delta = (\delta_1, \ldots , \delta_N)$ \\

$(X, \delta, \Lambda)$ & guided dynamical system with guiding sets
$\Lambda = (\Lambda_1, \ldots , \Lambda_N)$ \\

$\Phi_\delta$ & the semi-group of maps generated by $\delta$ \\

$\mathbf{id_X}$ & the identity map on $X$ \\

$OS(x)$ & the orbit set of a point $x$ \\

$\Lambda$-$OS(x)$ & the guided orbit set of a point $x$ \\

$\ell_p (\mathbb{R}^n)$ & $\mathbb{R}^n$ equipped with the norm
$\|\cdot\|_{\ell_p}$ \\

$\|\cdot\|_{\ell_p}$ & the norm $\|(x_1, \ldots, x_n) \|_{\ell_p}
= \left(\sum_{i=1}^n |x_i|^p \right)^{1/p}$ \\

$Df$ & the differential (Jacobian matrix) of a map $f$ \\

$C(X)$ & the space of all continuous functions on a topological
space $X$ \\

$C^k(M)$ & the space of $k$ times continuously differentiable
functions on $M$ \\

$L\left(\textbf{X},\textbf{Y}\right)$ & the space of bounded
linear operators from $\textbf{X}$ to $\textbf{Y}$ \\

$L\left(\textbf{X}\right)$ & the space of bounded linear operators
from $\textbf{X}$ into itself \\

$\textrm{Im}A$ & the image of a linear operator $A$ \\

$\textrm{Ker}A$ & the kernel $\{x | Ax = 0 \}$ of a linear
operator $A$ \\

$\textrm{ind}A$ & the index of a linear operator $A$ \\

$\textbf{I}$ & the identity operator on some function space \\

$\partial D$ & the boundary of a domain $D$ \\

$\overline{D}$ & the closure of a set $D$ \\

$\partial_x = \frac{\partial}{\partial x}$ & differentiation
with respect to the variable $x$ \\

$\partial = (\partial_x, \partial_y)$ & gradient
operator in the space $\mathbb{R}^2$ \\

$T_p (\Gamma)$ & the tangent space of the curve $\Gamma$ at the
point $p$ \\

$C_0^\infty (D)$ & the space of all infinitely differentiable
functions \\
 &  with compact support in $D$

\end{tabular}

\pagebreak

\chapter{Discrete guided dynamical systems with several
generators}\label{ch:dys}

In this chapter we shall present terminology and notation from the
theory of dynamical  systems essential for the formulation and
derivation of the results presented in later chapters. The
notation and terminology we use is not completely consistent with
the standard in this field. In particular, pay attention that we
use the term \emph{orbit-set} for what is usually called
\emph{orbit}, and the term \emph{orbit} will be reserved for a
more intuitive concept.

\section{Dynamical systems with several generators}\label{sec:definitions}
A \emph{dynamical system} is a pair $(X, \delta)$, where $X$ is a
metric space with a metric $d$ (usually compact) and $\delta =
(\delta_1 , \ldots , \delta_N)$ is a set of continuous maps
$\delta_i : X \rightarrow X$. The maps in $\delta$ generate (by
composition) a semigroup of maps $\Phi_\delta$ in the following
manner:
\begin{equation*}
\Phi_\delta^0 = \{\mathbf{id_X}\}
\end{equation*}
\begin{equation*}
\Phi_\delta^m = \{\sigma : X \rightarrow X | \exists  \sigma_1 ,
\ldots , \sigma_m \in \delta . \sigma = \sigma_1 \circ \cdots
\circ \sigma_m\}
\end{equation*}
and
\begin{equation*}
\Phi_\delta = \bigcup_{m=0}^\infty \Phi_\delta^m \,\,.
\end{equation*}
Given any $x_1 \in X$, an \emph{orbit} emanating from $x_1$ is a
sequence
\begin{equation*}
\mathcal{O} = (x_1, x_2, \ldots , x_n)
\end{equation*}
where for every $j = 2, \ldots, n$ there is some $i \in \{1,
\ldots , N\}$ such that
\begin{equation}\label{eq:orbit}
x_j = \delta_i (x_{j-1})
\end{equation}
We consider both finite and infinite orbits.

Given any $x \in X$, the \emph{orbit-set} of $x$ is the set
\begin{equation*}
OS(x) = \{\sigma(x) | \sigma \in \Phi_\delta\}
\end{equation*}
Equivalently, the orbit-set of a point $x$ may be defined as the
set of all $y$ for which there exists an orbit
\begin{displaymath}
\mathcal{O} = (x, \ldots, y)
\end{displaymath}
\begin{defins}
The following are basic notions relating to a dynamical system
$(X, \delta)$.
\begin{itemize}
  \item  $(X, \delta)$ is called \emph{minimal} if for
  all $x \in X$ it is true that $\overline{OS(x)} = X$.
  \item  A point $x_0 \in X$ is called an \emph{attractor} if there
  is a neighborhood $U$ of $x_0$ such that for any $x \in U$ there is
  an orbit emanating from $x$ and converging to $x_0$ \footnote{Some authors
  define an attractor as a \emph{set} having this property. Since we shall
  make no use of attractive sets which contain more than one point, we prefer
  to regard an attractor as a point.}.
  $x_0$ is called a \emph{global-attractor} if for
  any $x \in X$ there is an orbit emanating from $x$ and
  converging to $x_0$.
  \item  A point $x_0 \in X$ is called a \emph{weak attractor} if
  $x_0 \in \bigcap_{x \in X} \overline{OS(x)}$.
\end{itemize}
\end{defins}

\begin{rem} \emph{The term \emph{weak attractor} is not standard terminology
in dynamical systems. Nevertheless, this notion will prove to be
of key importance in the sequel, so the author took the right to
give this notion a name. Note that every global attractor is a
weak attractor. The reader with some experience in the general
theory of dynamical systems will note that a weak attractor is
nothing but a point lying in the intersection of the $\omega$ -
limit sets of all the points in $X$.}
\end{rem}
\begin{expl}
\emph{Let $X = [-1,1]\times[-1,1]$, put $p_1 = (1,1), p_2 =
(1,-1)$ , $p_3 = (-1,-1), p_4 = (-1,1)$. Define for $i = 1, 2, 3,
4$ the maps}
\begin{displaymath}
\delta_i: X \rightarrow X
\end{displaymath}
\emph{by}
\begin{displaymath}
\delta_i(x) = \frac{1}{2}(x + p_i)
\end{displaymath}
\emph{It follows from proposition \ref{prop:contracting} below
that $(X, \delta)$ is a minimal dynamical system. For any $i = 1,
2, 3, 4$, $p_i$ is an attractor - actually, a global attractor -
and these are the only attractors. On the other hand, as is the
case in any minimal system, any point in $X$ is a weak attractor.}
\end{expl}

It is useful to have at hand sufficient conditions for the
minimality of a dynamical system. The following proposition gives
one which will be useful later on.
\begin{prop} \label{prop:contracting}
Let $(X,d)$ be a compact, metric space, and let $\delta =
(\delta_1, \delta_2, \ldots , \delta_N )$ be a finite family of
functions $X \rightarrow X$ satisfying
\begin{equation}\label{eq:rangeCoverX}
\delta_1(X) \cup \delta_2(X) \cup \ldots \cup \delta_N(X) = X
\quad .
\end{equation}
If $\delta$ has the property that for all $ i = 1, \ldots , N$ and
all $x,y \in X$
\begin{equation}\label{eq:contracting}
x \neq y \ \Rightarrow \ d(\delta_i(x),\delta_i(y)) < d(x,y)
\end{equation}
then the dynamical system $(X,\delta)$ is minimal.
\end{prop}

\begin{proof}
Let us prove a lemma first.
\begin{lem}
For any $\epsilon > 0$ there exists a constant $0 \leq
c_{\epsilon} < 1$ such that for all $i = 1, \ldots , N$
\begin{equation*}
\forall x,y \in X \ . \ d(x,y) \geq \epsilon \Rightarrow
d(\delta_i(x),\delta_i(y)) \leq c_{\epsilon}d(x,y)
\end{equation*}
\end{lem}
\begin{proof}
Let there be given an $\epsilon
> 0$ and let $Y = X \times X$ with the product topology. For every $x \in
X$ let $B_\epsilon (x)$ denote the open ball around $x$ with
radius $\epsilon$. We define a compact subset $S \subseteq Y$ as
follows:
\begin{equation*}
S := Y \setminus \left[ \bigcup_{x \in X} \left( B_{\epsilon / 2}
(x) \times B_{\epsilon / 2} (x) \right) \right] \, \,.
\end{equation*}
For every $i = 1,2, \ldots , N$ define a function $g_i \colon S
\to \mathbb{R}$ by:
\begin{equation*}
g_i(x_1,x_2) = {d(\delta_i(x_1),\delta_i(x_2)) \over d(x_1,x_2)}
\end{equation*}
for all $(x_1,x_2) \in S$. For every $i$, $g_i$ is continuous, and
so $g_i$ attains a maximum  $c_{\epsilon , i}$. By
(\ref{eq:contracting}), \ $c_{\epsilon , i}<1$, for all $i$. Set
$c_\epsilon$ to be the maximum of these constants.

Now let $x,y$ be two points in $X$ s.t. $d(x,y) \geq \epsilon$.
Then we must have $(x,y) \in S$ so for every $i$
$$g_i(x,y) \leq c_\epsilon$$
and the lemma follows.
\end{proof}
Let us complete the proof of the proposition. Fix $x_0 \in X$ . To
prove the proposition we must show that for any $y$ in $X$ and
$\epsilon
> 0$ there is a $z \in OS(x_0)$ s.t. $d(z,y) \leq
\epsilon$. Fix some $y \in X$ and $\epsilon > 0$. Take some $n$
satisfying $c_{\epsilon}^n \cdot {\rm diam}(X) <\epsilon$, where
$c_{\epsilon}$ is the constant from the lemma. The lemma tells us
that for all $\sigma \in \Phi_\delta^n$ and all $x_1, x_2 \in X$
\begin{displaymath}
d(\sigma(x_1),\sigma(x_2)) \leq \epsilon
\end{displaymath}
and thus for all $\sigma \in \Phi_\delta^n$ :
\begin{equation} \label{eq:diam}
{\rm diam}(\sigma(X)) \leq \epsilon
\end{equation}

But note that by virtue of (\ref{eq:rangeCoverX}),
\begin{displaymath}
\bigcup_{f \in \Phi_\delta^n} f(X) = X
\end{displaymath}
so that there is an $f \in \Phi_\delta^n$ s.t. $y \in f(X)$. Now
by (\ref{eq:diam}) it follows that for all $x$ it is true that
$d(f(x),y) \leq \epsilon$ so we can choose $z = f(x_0)$ and the
proof is complete.
\end{proof}

\begin{defins}
Let $(X, \delta)$ be a dynamical system.
\begin{itemize}
  \item A set $Y \subseteq X$ is called \emph{$\delta$-invariant}
  if $\delta_i(y) \in Y$ for all $i = 1, \ldots , N$ and
  $y \in Y$.
  \item If $Y \neq \emptyset$ is a closed, $\delta$-invariant subset of $X$ then
  $\delta$ naturally induces on $Y$ a dynamical
  system $(Y, \tilde{\delta})$, where
  \begin{displaymath}
  \tilde{\delta} = (\delta_1 |_Y , \ldots , \delta_N |_Y )
  \end{displaymath}
  $(Y, \tilde{\delta})$ is called a \emph{subsystem} of $(X, \delta)$. Because there is no
  chance of ambiguity, we shall denote this dynamical system simply by
  $(Y, \delta)$.
\end{itemize}
\end{defins}

\section{Guided dynamical systems}\label{sec:guideddys}

Usually, in the study of dynamical systems, one is interested in
the behavior of points under the action of $\Phi_\delta$, that is,
``how a point moves" under iterations of maps in $\Phi_\delta$.
Such movement may be described by the class of all orbits of
point. But in certain applications of dynamical systems it is most
profitable to ignore certain, ``illegal'', orbits and to
concentrate on a subclass of the orbits. These ideas were
introduced by Paneah in \cite{Pa03NC}, \cite{Pa03SFE} and
\cite{Pa04DAIG}, and will be developed below.

\begin{defin}
A \emph{guided dynamical system} is a dynamical system
\newline $\left(X, ( \delta_1, \ldots , \delta_N ) \right)$ together with a
system $\Lambda = (\Lambda_1 , \ldots , \Lambda_N )$ of $N$ closed
subsets of $X$.
\end{defin}
The sets $\Lambda_i$ are called \emph{guiding sets}. It will be
always assumed that $\bigcap_{i=1}^N \Lambda_i = \emptyset$. We
shall also denote at some times the set $\bigcup_{i=1}^N
\Lambda_i$ by $\Lambda$. This will never cause any confusion.

\begin{defin}
An orbit is called a \emph{$\Lambda$-proper orbit} , or, for
short, a \emph{$\Lambda$-orbit} , if in (\ref{eq:orbit}) $\delta_i
\neq \delta_k$ if $x_{j-1} \in \Lambda_k$.
\end{defin}
When studying a guided dynamical system we restrict our attention
to $\Lambda$-proper orbits. One can think of a guided dynamical
system as a dynamical system with several generators in which
there are points that one can leave using only a subset of
$\delta$. A different point of view is to consider $\delta_i$ as a
function with a domain of definition $X \setminus \Lambda_i$. For
true motivation for this concept the reader must wait until
chapters 2 and 4.
\begin{rem}
\emph{When dealing with a dynamical system with only two
generators, $\Lambda_1$ is the set of points which we \emph{must}
leave using $\delta_2$, and vice versa. So one may equivalently
define $\mathcal{T}_1 = \Lambda_2$ and $\mathcal{T}_2 = \Lambda_1$
to be the guiding sets, that is, to associate the guiding the set
with the map which we must use on it. Actually, this is the
original notation used by Paneah.}
\end{rem}

\begin{metadefin}
Let $(X, \delta, \Lambda)$ be a guided dynamical system, and let
$\clubsuit$ be some concept relating to the dynamical system $(X,
\delta)$ which may be defined by means of the orbits in $(X,
\delta)$. Then $\Lambda$-$\clubsuit$ is the concept relating to
the guided dynamical system $(X, \delta, \Lambda)$ which is
defined precisely as $\clubsuit$ with the difference that the
phrase ``orbit'' is replaced by the phrase ``$\Lambda$-proper
orbit''.
\end{metadefin}
For example, the \emph{$\Lambda$-orbit set} of a point $x$,
denoted $\Lambda$-$OS(x)$, is the set of points $y$ for which
there exists a \emph{$\Lambda$-proper} orbit $\mathcal{O}$
\begin{displaymath}
\mathcal{O} = (x, \ldots, y)
\end{displaymath}
We may similarly define a $\Lambda$-attractor, a $\Lambda$-minimal
dynamical system, etc.

\begin{expl}\label{expl:circle}
\emph{Let $X = \mathbb{S}^1$, the unit circle in the complex
plane, and let}
\begin{eqnarray*}
\delta_1 (z) & = & e^{i 2 \pi\theta_1} z \\ \delta_2 (z) & = &
e^{i 2 \pi\theta_2} z
\end{eqnarray*}
\emph{Define $\Lambda_1 = \{1, -1\}$ and $\Lambda_2 = \{i, -i\}$.
By the well known theorem of Kronecker and Weyl
(\cite{MathWorld}), the dynamical system $(X, \delta)$, (when
viewed as an unguided dynamical system), is minimal if and only if
at least one of $\theta_1, \theta_2$ is irrational. Does this
remain true when $(X, \delta, \Lambda)$ is viewed as a
\emph{guided} dynamical system? Let us show that the answer to
this is almost yes. To be precise, we shall show that $(X, \delta,
\Lambda)$ is $\Lambda$-minimal if and only if at least one of
$\theta_1, \theta_2$, say $\theta_1$, is irrational and the other
one, say $\theta_2$, is not an integer multiple of $\frac{1}{2}$.}

\emph{The ``only if'' part is clear. Now assume, without loss of
generality, that $\theta_1 \notin \mathbb{Q}$. We also assume that
$\theta_2 \in \mathbb{Q}$, as the proof in the case $\theta_2
\notin \mathbb{Q}$ is similar. Let $z_1$ be a point on the circle.
We have to prove that} $\overline{\Lambda-OS(z_1)} =
\mathbb{S}^1$.

\emph{Consider the maximal $\Lambda$-proper orbit of the type}
\begin{displaymath}
\mathcal{O} = \left(z_1, \delta_1 (z_1), \delta_1 (\delta_1
((z_1)), \ldots \right)
\end{displaymath}
\emph{Denote this maximal orbit by $\tilde{\mathcal{O}}$. Note
that $\tilde{\mathcal{O}}$ is at least one point long. There are
only two possibilities:}
\begin{enumerate}
  \item $\tilde{\mathcal{O}}$ is infinite (this happens when
  $\tilde{\mathcal{O}}$ never intersects $\Lambda_1$). In this
  case, by the Kronecker-Weyl theorem, $\tilde{\mathcal{O}}$ is
  dense in $\mathbb{S}^1$, so $\Lambda\verb"-"OS(z_1)$ is, too.
  \item $\tilde{\mathcal{O}}$ is finite. This means that for some
  $m \in \mathbb{N}$, $\delta_1^m (z_1) \in \Lambda_1$ \footnote{By $\delta_1^m$
  we mean the $m$th iterate of $\delta_1$.}, but also $\delta_1^m (z_1) \in
  \Lambda\verb"-"OS(z_1)$. Now, $\theta_2$ is not an integer
  multiple of $\frac{1}{2}$, so $\delta_2(\delta_1^m(z_1)) \notin
  \Lambda_1$, and
  \begin{displaymath}
    (\delta_2(\delta_1^m (z_1)), \delta_1 ( \delta_2(\delta_1^m
    (z_1))), \delta_1(\delta_1 ( \delta_2(\delta_1^m
    (z_1)))), \ldots )
  \end{displaymath}
    is now an infinite orbit that doesn't intersect
    $\Lambda_1$, therefore it is dense in $\mathbb{S}^1$. Because $\delta_1^m (z_1) \in
    \Lambda-OS(z_1)$
    this implies that the orbit set of $z_1$ is dense in
    $\mathbb{S}^1$.
\end{enumerate}
\emph{Examining the above proof one sees that even if $\theta_2$
is equal to $\frac{1}{2}$, the points $z=1$ and $z=-1$ are
$\Lambda$-weak attractors if $\theta_1$ is irrational.}
\end{expl}

The next concept we shall introduce turns out to be crucial for
stating necessary and sufficient conditions for unique solvability
of functional equations and boundary value problems , so we shall
be explicit when defining it.
\begin{defin}
A set $Y \subseteq X$ is called \emph{  $(\Lambda,
\delta)$-invariant} if
\begin{displaymath}
\forall y \in Y . \forall i . y \notin \Lambda_i \Rightarrow
\delta_i(y) \in Y
\end{displaymath}
In words, any $\Lambda$-orbit that begins in $Y$ also ends there.
\end{defin}

It is a well known fact in the theory of dynamical systems that
any compact dynamical system \footnote{By  \emph{compact dynamical
system} we mean a dynamical system $(X, \delta)$ where $X$ is
compact.} $(X, \delta)$ has a closed subsystem $(A, \delta)$ that
is minimal (see \cite{HaKatok}). It is interesting to note that
with some care this result carries over to guided dynamical
systems as well.

\begin{lem}\label{lem:invariantclosure}
Let $(X, \delta)$ be a dynamical system, and let $Y$ be a
$(\Lambda, \delta)$-invariant subset of $X$. Then $\overline{Y}$
is also $(\Lambda, \delta)$-invariant.
\end{lem}
\begin{proof}
Let $y \in \overline{Y}$, and assume that $I\subseteq \{1, \ldots
, N\}$ is the set of indices $i$ for which $y \notin \Lambda_i$.
we have to show that
\begin{displaymath}
\forall i \in I . \delta_i(y) \in \overline{Y}
\end{displaymath}
Fix $i \in I$. There is a sequence $(y_n)_{n=1}^\infty$ of points
in $Y$ such that $y_n \rightarrow y$. Since $\Lambda_i$ is closed,
for sufficiently large $n $, $y_n \notin \Lambda_i$. Since $Y$ is
$(\Lambda, \delta)$-invariant, for these $n$ we have
$\delta_i(y_n) \in Y$. By continuity of $\delta_i$, $\delta_i(y_n)
\rightarrow \delta_i(y)$, so $\delta_i(y) \in \overline{Y}$, as
required. Now since $i$ was an arbitrary element of $I$, the proof
is complete.
\end{proof}

\begin{thm}\label{thm:existminimal}
Every compact guided dynamical system $(X, \delta, \Lambda)$ has a
closed, $\Lambda$-minimal, $(\Lambda, \delta)$-invariant
subsystem.
\end{thm}
\begin{proof}
Let $(X, \delta, \Lambda)$ be a compact guided dynamical system.
Denote by $\mathcal{M}$ the collection consisting of all closed,
non-empty, $(\Lambda, \delta)$-invariant subsets of $X$.
$\mathcal{M}$ is not empty, because $X \in \mathcal{M}$. We shall
use Zorn's lemma to prove that $\mathcal{M}$ has a
minimal\footnote{Here, of course, we are using the word
\emph{minimal} in the usual sense, that is, minimal with respect
to inclusion.} element.

Assume that
\begin{displaymath}
\left\{A_\alpha \right\}_\alpha
\end{displaymath}
is a chain in $\mathcal{M}$. A lower bound for this chain is given
by
\begin{displaymath}
B \triangleq \bigcap_\alpha A_\alpha
\end{displaymath}
Indeed, let us prove that $B \in \mathcal{M}$. Obviously, $B$ is
closed. Also, $B \neq \emptyset$, because if it is empty then, $X$
being compact, there must $A_{\alpha_1}, \ldots , A_{\alpha_M}$
such that
\begin{displaymath}
\bigcap_{k=1}^M A_{\alpha_k} = \emptyset
\end{displaymath}
But the above intersection is decreasing and thus equals one of
the $A_\alpha$'s, contradicting the assumption that for all
$\alpha$, $A_\alpha \neq \emptyset$. Finally, $B$ is $(\Lambda,
\delta)$-invariant. Indeed, let $b \in B$, and assume that $I
\subseteq \{1, \ldots , N \}$ is the set of indices $i$ such that
$b \notin \Lambda_i$. For all $\alpha$ and all $i \in I$, $b \in
A_\alpha$ and $b \notin \Lambda_i$. By the $(\Lambda,
\delta)$-invariance of $A_\alpha$ we have that $\delta_i(b) \in
A_\alpha$. This is true for all $\alpha$, so $\delta_i(b) \in
\bigcap_\alpha A_\alpha = B$. Since this is true for all $i \in
I$, $B$ is $(\Lambda, \delta)$-invariant, and thus is in
$\mathcal{M}$.

Now Zorn's lemma guaranties the existence of a closed, non-empty,
$(\Lambda, \delta)$-invariant $A \subseteq X$. It is left to show
that $A$ is $\Lambda$-minimal.

Take any $x \in A$. $\Lambda$-$OS(x)$ is definitely $(\Lambda,
\delta)$-invariant. By the previous lemma, so is
$\overline{\Lambda-OS(x)}$. As $A$ is invariant,
$\overline{\Lambda-OS(x)} \subseteq A$. By the minimality of $A$,
\begin{displaymath}
\overline{\Lambda-OS(x)} = A
\end{displaymath}
and, since $x$ was arbitrary in this discussion, this means that
$(A, \delta, \left(\Lambda_1 \cap A, \ldots, \Lambda_N \cap A
\right))$ is a minimal dynamical system.
\end{proof}

\begin{prop}\label{prop:miniffno}
A guided dynamical system  $(X, \delta, \Lambda)$ is
$\Lambda$-minimal if and only if it has no $\Lambda$-subsystem
other than itself.
\end{prop}
\begin{proof}
Taking into account  \ref{lem:invariantclosure} and the fact that
a subsystem is nothing but a closed, non-empty, invariant subset,
the assertion is clear.
\end{proof}

\section{Isomorphism of guided dynamical systems}\label{sec:iso}

For every abstract mathematical structure it is always useful to
define the maps between two instances of the same type of
structure that preserve the essential features of that structure.
In the standard theory of dynamical systems, there are the
important concepts of a \emph{factor} and an \emph{isomorphism} of
dynamical system. More details are to be found in \cite{HaKatok}.
We shall restrict our attention only to isomorphism of two
(guided) dynamical systems, as this term will be very useful later
on.
\begin{defin} \label{defin:iso}
Two dynamical systems $(X, (\delta_1, \ldots , \delta_N ))$
and $(Y, (\gamma_1, \ldots , \gamma_N ))$ are said to be
\emph{isomorphic} if there exists a homeomorphism $\varphi : X
\rightarrow Y$ satisfying
\begin{equation*}
\varphi \circ \delta_i \circ \varphi^{-1} = \gamma_i \quad {\rm
for} \quad i = 1, \ldots, N
\end{equation*}
$\varphi$ is called an \emph{isomorphism} , of the dynamical
systems $(X, \delta)$ and $(Y, \gamma)$.
\end{defin}
Loosely speaking, isomorphic dynamical systems exhibit the same
dynamical behavior. For instance, $x \in X$ is an attractor if and
only if $\varphi(x)$ is an attractor in $(Y, \gamma)$, and $(X,
\delta)$ is minimal if and only if $(Y, \gamma)$ is minimal, and
so on.
\begin{defin}
Two guided dynamical systems $(X, \delta, (\Lambda_1, \ldots,
\Lambda_N))$ and \newline $(Y, \gamma, (\Omega_1, \ldots,
\Omega_N))$ are said to be \emph{isomorphic} if $(X, \delta)$ and
$(Y, \gamma)$ are isomorphic as dynamical systems and $\varphi$
from definition \ref{defin:iso} maps each $\Lambda_i$ onto
$\Omega_i$.
\end{defin}

For completeness of this exposition, let us prove two results
regarding isomorphic guided dynamical systems.

\begin{lem}\label{lem:propiffprop}
Let $(X, \delta, \Lambda)$ and $(Y, \gamma, \Omega)$ be two guided
dynamical systems. If $\varphi : X \rightarrow Y$ is an
isomorphism of guided dynamical systems then the orbit
\begin{equation*}
\mathcal{O} = (x_1, x_2, \ldots , x_n)
\end{equation*}
is $\Lambda$-proper if and only if
\begin{equation*}
\tilde{\mathcal{O}} = (\varphi(x_1), \varphi(x_2), \ldots ,
\varphi(x_n))
\end{equation*}
is $\Omega$-proper.
\end{lem}
\begin{proof}
Note that $(x_1, x_2, \ldots , x_n)$ is $\Lambda$-proper if and
only if $(x_1, x_2)$, $(x_2, x_3)$, \ldots , $(x_{n-1}, x_n)$ are
all $\Lambda$-proper. So we may assume that $\mathcal{O} = (x_1,
x_2)$. Also, by the symmetry of the relation ``isomorphic", it
suffices to show that $\Lambda$-properness of $\mathcal{O}$
implies $\Omega$-properness of $\tilde{\mathcal{O}}$.

Assume then that $\mathcal{O} = (x_1, x_2)$ is $\Lambda$-proper.
We must have $x_2 = \delta_i (x_1)$, for some $i \in \{1, \ldots ,
N\}$, and $x_1 \notin \Lambda_i$. Because $\varphi$ is an
isomorphism of dynamical systems
\begin{displaymath}
\varphi(x_2) = \varphi(\delta_i (x_1)) = \gamma_i (\varphi(x_1))
\end{displaymath}
and this shows that $\tilde{\mathcal{O}} = (\varphi(x_1),
\varphi(x_2))$ is an orbit. To see that it is $\Omega$-proper, we
just note that as $\varphi$ is a 1-1 function that maps guiding
sets onto guiding sets and $x_1 \notin \Lambda_i$, $\varphi(x_1)$
cannot be in $\Omega_i$.
\end{proof}

\begin{thm}\label{thm:isomorphic}
Let $(X, \delta, \Lambda)$ and $(Y, \gamma, \Omega)$ be two
isomorphic guided dynamical systems.
\begin{enumerate}
  \item $x_0$ is a $\Lambda$-weak attractor in $(X, \delta, \Lambda)$
  if and only if $\varphi(x_0)$ is an $\Omega$-weak attractor in $(Y, \gamma, \Omega)$.
  \item $(X, \delta, \Lambda)$ is $\Lambda$-minimal if and only if
    $(Y, \gamma, \Omega)$ is $\Omega$-minimal.
\end{enumerate}
\end{thm}
\begin{proof}
Assume that $x_0$ is a $\Lambda$-weak attractor in $(X, \delta,
\Lambda)$. Let $y \in Y$. We must show that $\varphi(x_0) \in
\overline{\Omega-OS(y)}$. Choose any $\epsilon > 0$ and define $x
= \varphi^{-1} (y)$. By the continuity of $\varphi$, there is a
$\mu
> 0$ such that $d_X (z,x_0)< \mu$ implies $d_Y (\varphi(z),
\varphi(x_0))< \epsilon$ \footnote{Here $d_X$ and $d_Y$ denote the
metrics on $X$ and $Y$, respectively.} for all $z \in X$. There is
a $\Lambda$-proper orbit in $X$
\begin{displaymath}
\mathcal{O}=(x, x_1, \ldots, x_n)
\end{displaymath}
such that $d_X (x_n, x_0) < \mu$. But then by the lemma
\begin{displaymath}
\tilde{\mathcal{O}}=(\varphi(x) = y, \varphi(x_1), \ldots,
\varphi(x_n))
\end{displaymath}
is $\Omega$-proper and $d_Y (\varphi(x_n), \varphi(x_0))
<\epsilon$. Since this argument is valid for any $\epsilon > 0$,
we have that $\varphi(x_0) \in \overline{\Omega-OS(y)}$. We have
established the ``only if" half of the first part of the theorem.
The ``if" part follows by interchanging the roles of $(X, \delta,
\Lambda)$ and $(Y, \gamma, \Omega)$. Finally, the second part of
the theorem clearly follows from the first.
\end{proof}

\chapter{Some results in functional equations}\label{ch:FuncEqs}

In the following sections we will show how the notions and results
of chapter 1 are applied in the field of functional equations. We
shall not attempt to explain what a functional equation is, the
history of functional equations, and so forth. Such information
may be found in the fundamental works of two of the leading
specialist in this field in the 20th century: Janos Acz\'{e}l
(\cite{Aczel1}, \cite{Aczel2}) and Marek Kuczma (\cite{Kuczma1},
\cite{Kuczma2}).

\section{The Maximum principle for functional
equations}\label{sec:Maxp} Maximum principles appeared in analysis
long ago. In the theory of functions of a complex variable, the
maximum modulus principle for analytical functions helps to
establish further results - e.g. Schwartz's lemma.  
In partial differential equations they serve as a tool for proving
uniqueness theorems, approximating solutions, etc.
A maximum principle in the field of functional
equations appeared for the first time only a few years ago. In
2003 Paneah showed in \cite{Pa03SFE} and \cite{Pa03NC} that under
certain assumptions, if a function $F$ satisfies
\begin{displaymath}
F(t) - a_1(t) F(\delta_1(t)) - a_2(t) F(\delta_2(t)) = 0 \hspace{5
mm} , \hspace{5 mm} t \in [-1,1]
\end{displaymath}
then $F$ attains its maximum and minimum values on the boundary of
$[-1,1]$. This theorem proved useful for applications in integral
geometry, partial differential equations and, of course, in
functional equations (\cite{Pa04DMFE} and \cite{Pa04DAIG}).

The purpose of this section is to extend Paneah's maximum
principle as far as we can in order to prove a uniqueness theorem
for a conditional cauchy equation in $\mathbb{R}^n$. Throughout
this section $(X,\delta, \Lambda)$ will be a guided dynamical
system.
\subsection{The maximum principle}\label{subsec:maxp}
To begin with, let us recall the notion of a semi-continuous
function.
\begin{defin}
Let $X$ be a metric space and $x_0 \in X$. A real valued function
$f: X \rightarrow \mathbb{R}$ is said to be \emph{upper
semi-continuous at $x_0$} if
\begin{displaymath}
\limsup_{x \rightarrow x_0} f(x) \leq f(x_0) \quad .
\end{displaymath}
$f$ is a said to be \emph{upper semi-continuous} if it is upper
semi-continuous at any point $x \in X$. A real valued function $f$
is called \emph{lower semi-continuous (at a point $x_0$)} if $x
\mapsto -f(x)$ is upper semi-continuous (at the point $x_0$).
\end{defin}
\begin{lem} \label{lem:maxp}
Let $f : X \rightarrow \mathbb{R}$ be an upper (lower)
semi-continuous function that satisfies the following functional
equation:
\begin{equation}\label{eq:funceq1}
f(x) - \sum_{i=1}^N a_i (x) \cdot f(\delta_i (x)) = 0 \verb" , " x
\in X
\end{equation}
where $a_i : X \rightarrow \mathbb{R}$ satisfy :
\begin{eqnarray}\label{eq:coef}
  \forall i . \forall x . a_i(x) \geq 0 \\ \label{eq:coef1}
  \forall i . \forall x \notin \Lambda_i . a_i (x) > 0 \\ \label{eq:coef2}
  \forall x . \sum_{i=1}^N a_i(x) = 1 .
\end{eqnarray}
Then if $f$ attains its maximum (minimum) at some point $y_0 \in
X$, then it attains its maximum (minimum) at any point $x \in
\overline{\Lambda-OS(y_0)}$.
\end{lem}
\begin{proof}
Put $M = f(y_0) = \max f$, and let $I \subseteq \{1, \ldots , N\}$
be a subset of indices $i$ such that $y_0 \notin \Lambda_i$. Then
there are numbers $\epsilon_1, \ldots, \epsilon_N \geq 0$ such
that
\begin{equation*}
f(\delta_i (y_0)) = M - \epsilon_i , \quad i=1, \ldots, N \,.
\end{equation*}
Combining these relations with (\ref{eq:funceq1}) and using
(\ref{eq:coef}), (\ref{eq:coef1}), (\ref{eq:coef2}) results in
\begin{equation*}
\sum_{i \in I} a_i(y_0)\cdot \epsilon_i = 0 \,.
\end{equation*}
Thus $\epsilon_i = 0$ and so $f(\delta_i(y_0)) = M$  for all $i
\in I$. Now by induction, for any point $x \in \Lambda-OS(y_0)$ we
have $f(x) = M$. If $x \in \overline{\Lambda-OS(y_0)}$ then there
is a sub-sequence $x_n \rightarrow x$ from $\Lambda-OS(y_0)$, and
since $f$ is upper semi-continuous, we have

$$ M = \lim_{n \rightarrow \infty} M = \lim_{n \rightarrow \infty} f(x_n) \leq f(x) $$
so $f(x) = M$, which was to be proved.
\end{proof}

\begin{cor}
Let $(X, \delta, \Lambda)$ be a compact, $\Lambda$-minimal
dynamical system. Assume that $f : X \rightarrow \mathbb{R}$ is an
upper semi-continuous function that satisfies (\ref{eq:funceq1})
where the coefficients $a_i$ satisfy (\ref{eq:coef}),
(\ref{eq:coef1}), (\ref{eq:coef2}). Then $f$ is constant.
\end{cor}
\begin{proof}
Being an upper semi-continuous function on a compact space, $f$
attains a maximum $M = \max_X f$ at some point $y_0 \in X$. The
$\Lambda$-minimality of $(X, \delta, \Lambda)$ gives us
$\overline{\Lambda-OS(y_0)} = X$. Using the lemma we assert that
$f \equiv M$.
\end{proof}

We now proceed to prove a lemma which will be useful when proving
the main result of this section.
\begin{lem}\label{lem:lemmaxp}
Assume that $(X, \delta, \Lambda)$ is a compact guided dynamical
system having a $\Lambda$-weak attractor $x_0 \in X$. Assume that
$f : X \rightarrow \mathbb{R}$ is a continuous solution of
equation (\ref{eq:funceq1}) where all the coefficients $a_i$
satisfy relations (\ref{eq:coef}), (\ref{eq:coef1}),
(\ref{eq:coef2}). Then the function $f$ is constant.
\end{lem}
\begin{proof}
As the function $f$ is continuous, there are points $y_0, y_1 \in
X$ for which $f(y_0) = \min_X f$ and $f(y_1) = \max_X f$. Being a
$\Lambda$-weak attractor, the point $x_0$ belongs to both sets
\begin{displaymath}
\overline{\Lambda-OS(y_0)} \,\,\,{\rm and}\,\,\,
\overline{\Lambda-OS(y_1)} \,\,.
\end{displaymath}
Being continuous, the function f is simultaneously upper and lower
semi-continuous, and hence, by lemma \ref{lem:maxp}, it takes its
maximum and minimal values at $x_0$. It follows that $f \equiv
f(x_0) = const $, and this completes the proof of the lemma.
\end{proof}

\begin{expl}
\emph{Let $\mathbb{S}^1$ be the unit circle in the complex plane.
Consider the functional equation}
\begin{equation}\label{eq:exmplFE}
f(z) = \sin^2 (\arg z)f(e^{i \tau_1} \cdot z) + \cos^2 (\arg
z)f(e^{i \tau_2} \cdot z) \hspace{5 mm} , \hspace{5 mm} z \in
\mathbb{S}^1
\end{equation}
\emph{where $\tau_1, \tau_2 \in \mathbb{R}$ are fixed constants.
We claim that equation (\ref{eq:exmplFE}) has a non-constant
continuous solution if and only if both numbers $\tau_1 / 2\pi$
and $\tau_2 / 2\pi$ are rational. Indeed, if $\tau_1 / 2\pi,
\tau_2 / 2\pi \in \mathbb{Q}$ then we may write}
\begin{eqnarray*}
\tau_1 = {2 \pi k_1}/n \\
\tau_2 = {2 \pi k_2}/n
\end{eqnarray*}
\emph{with $k_1 , k_2 , n \in \mathbb{Z}$.Then for an arbitrary
continuous $(\frac{2\pi}{n})-$periodic function $g : \mathbb{R}
\rightarrow \mathbb{R}$ the function}
\begin{displaymath}
f(z) = g(\arg z)
\end{displaymath}
\emph{is a continuous solution of (\ref{eq:exmplFE}).}

\emph{On the other hand, if, for example, $\tau_1 / 2\pi \notin
\mathbb{Q}$, then as we have shown in example \ref{expl:circle},
the guided dynamical system on the circle generated by the
functions  $\delta_1 (z) = e^{i\tau_1} \cdot z$ and $\delta_2 (z)
= e^{i\tau_2} \cdot z$ and guided by $\Lambda_1 = \{z | \sin^2
(\arg z) = 0 \} = \{1, -1 \}$ and $\Lambda_2 = \{z | \cos^2 (\arg
z) = 0 \} = \{i, -i \}$ has the point $z=1$ is a $\Lambda$-weak
attractor. Thus by lemma \ref{lem:lemmaxp} equation
(\ref{eq:exmplFE}) has no non-constant continuous solutions, the
requirements on the coefficients being clearly fulfilled.}
\end{expl}
\begin{thm}\label{thm:maxp}
Assume that $(X, \delta, \Lambda)$ is compact guided dynamical
system that has a $\Lambda$-weak attractor $x_0 \in X$. Assume
that a function $F : X \rightarrow \mathbb{R}^n$ is a continuous
solution of the equation
\begin{equation}\label{eq:funceqvec}
F(x) - \sum_{i=1}^N A_i (x) \cdot F(\delta_i (x)) = 0  \verb" , "
x \in X \,\, .
\end{equation}
The coefficients $A_i : X \rightarrow \mathbb{R}^{n \times n}$ are
assumed to be lower triangular matrices with non-negative entries
on the diagonal for all $x \in X$ and satisfy:
\begin{eqnarray}\label{eq:coefvec}
\forall i . \forall x \notin \Lambda_i \det(A_i(x))>0\label{eq:coefvec1}\\
\forall x . \sum_{i=1}^N A_i(x) = I \label{eq:coefvec2} \,\,.
\end{eqnarray}
Then $F$ is constant.
\end{thm}
\begin{proof}
We write $F(x) = (f_1 (x), \ldots, f_n (x))$, with $f_k : X
\rightarrow \mathbb{R}$ continuous for every $k$. We also write
$A_i^{k,m}$ for the entry in the $k$th row and $m$th column in the
matrix $A_i (x)$. Equation (\ref{eq:funceqvec}) may now be written
as a system of $n$ functional equations:
\begin{equation}\label{eq:system}
f_k (x) - \sum_{m=1}^k A_1^{k,m} f_m (\delta_1 (x)) - \ldots -
\sum_{m=1}^k A_N^{k,m} f_m (\delta_N (x)) = 0
\end{equation}
for $k = 1, \ldots , n$. The first equation is
\begin{equation}\label{eq:systemfirst}
f_1(x) - A_1^{1,1} f_1 (\delta_1 (x)) - \ldots - A_N^{1,1} f_1
(\delta_N (x)) = 0
\end{equation}
and, using (\ref{eq:coefvec1}) and (\ref{eq:coefvec2}) the lemma
tells us that $f_1 \equiv c_1$.

Assume that $f_1 \equiv c_1, \ldots , f_k \equiv c_k$. Let us show
that $f_{k+1} \equiv c_{k+1}$. Indeed, for $k+1$ we may rewrite
(\ref{eq:system}) as
\begin{equation}\label{eq:systeminduction}
f_{k+1} (x) - \sum_{i=1}^N A_i^{k+1,k+1} f_{k+1} (\delta_i (x)) -
\sum_{m=1}^k \sum_{i=1}^N A_i^{k+1,m} c_m = 0 .
\end{equation}
But (\ref{eq:coefvec2}) means that $\sum_{i=1}^N A_i^{k+1,m} c_m =
0$ for $m < k+1$ so (\ref{eq:systeminduction}) reduces to
(\ref{eq:funceq1}) and again the lemma ensures that $f_{k+1}
\equiv c_{k+1}$. This completes the proof of the theorem.
\end{proof}

Assume now that for every $x$, $\{A_1(x), \ldots, A_N(x)\}$ is a
commuting family of matrices with only (real) positive
eigenvalues. A basic result from linear algebra says that for
every $x$ there is an invertible matrix $P_x \in \mathbb{R}^{n
\times n}$ such that $T_i (x) = P^{-1}_x A_i(x) P_x$ is lower
triangular for all $i$ \footnote{\cite{HoffKun}}. In some (very
rare, unfortunately) cases, $P_x$ can be chosen to be constant
throughout $X$, that is, $P_x = P$. Assume that this is the case.
If $F$ satisfies (\ref{eq:funceqvec}) we may equivalently write
that equation as:
\begin{equation*}
P P^{-1} \cdot F(x) - \sum_{i=1}^N P T_i (x) P^{-1} \cdot
F(\delta_i (x)) = 0  \verb" , "  x \in X
\end{equation*}
or
\begin{equation}\label{eq:Pfunceqvec}
P [P^{-1} \cdot F(x) - \sum_{i=1}^N T_i (x) P^{-1} \cdot
F(\delta_i (x))] = 0  \verb" , "  x \in X .
\end{equation}
We define a new function $G(x) = P^{-1} \cdot F(x)$. Because $P$
is invertible, (\ref{eq:Pfunceqvec}) can be re-written as
\begin{equation*}
G(x) - \sum_{i=1}^N T_i (x) \cdot G(\delta_i (x)) = 0 \verb" , " x
\in X
\end{equation*}
and this is exactly the situation of theorem \ref{thm:maxp}. We
record this result as
\begin{cor}\label{cor:matrixCTFE}
Let the assumptions of theorem \ref{thm:maxp} hold with the single
change that now $A_1(x), \ldots, A_N(x)$ form a commuting family
of matrices with only (real) positive eigenvalues for which there
exists a constant triangulating matrix (that is good for all $x$).
Then $F$ is constant.
\end{cor}

In the above discussion we assumed the existence of a matrix $P$
that triangulates $A_i (x)$ for all $x \in X, i=1, \ldots, N$.
When can we be sure that such a matrix exists? Trivially, when
$A_i(x)$ are already triangular. Also, if $A_i (x) =
\varphi_i(x)B_i$, where $B_i$ is a constant matrix for all $i$,
and $\varphi_i$ is some real valued function, we can find a
constant matrix $P$ that does the job. However, the latter class
of matrix-functions will never arise non-trivially in our
applications.

\subsection{An application to Cauchy type functional equations}\label{subsec:CauchyRn}
We now use the results of the previous section to find the $C^1$
\footnote{Given a compact subset $K$ of $\mathbb{R}^n$, we denote
by $C^1(K)$ the space of all functions on $K$ that have
continuously differentiable extensions to every neighborhood $U$
of $K$.} solutions to certain functional equations of the type
\begin{equation}\label{eq:cauchyveceq}
f(x) - f(a_1(x)) - f(a_2(x)) = 0 \verb"  ,  " x \in \mathbf{K} .
\end{equation}
Following Paneah we shall call equations of the above type
\emph{Cauchy type functional equations}.

\begin{thm}\label{thm:cauchyveceq}
Let $\mathbf{K}$ be a compact, connected subset of $\mathbb{R}^n$.
Let $a_1, a_2 : \mathbf{K} \rightarrow \mathbf{K} $ be $C^1$ maps
that generate a dynamical system in $\mathbf{K}$ with a weak
attractor and satisfy
\begin{equation*}
\forall x \in \mathbf{K} . a_1(x) + a_2(x) = x
\end{equation*}
Assume that the differentials $A_1(x)$ and $A_2(x)$ of $a_1(x)$
and $a_2(x)$ have only (real) positive eigenvalues. Assume also
that there exists an invertible matrix $P$ such that for any point
$x \in X$ there exists two lower triangular matrices $T_1 (x)$ and
$T_2 (x)$ such that
\begin{displaymath}
T_i (x) = P^{-1} A_i(x) P \, , \,\,\, i=1,2 \,\,.
\end{displaymath}
If $f \in C^1 (\mathbf{K},\mathbb{R})$ is a solution of
(\ref{eq:cauchyveceq}), then there exists a vector $c \in
\mathbb{R}^n$ such that
\begin{displaymath}
 f(x) = c \cdot x \quad .
\end{displaymath}
\end{thm}
\begin{proof}
Assume that $f$ is a solution of (\ref{eq:cauchyveceq}). Denote by
$\nabla f(x) $ , $A_1 (x)$ and $A_2(x)$ the differentials of $f$
$a_1$ and $a_2$ , respectively, at the point $x$. Put $g(x) =
(\nabla f(x))^T$ and $B_i (x) = (A_i (x))^T$. Then differentiating
(\ref{eq:cauchyveceq}) we obtain
\begin{equation*}
g(x) - B_1 (x) \cdot g(a_1(x)) - B_2 (x) \cdot g(a_2(x)) = 0 \quad
.
\end{equation*}
Note that $B_1 = I_{n \times n} - B_2$, so that these matrices
commute, and they also have exactly the same eigenvalues as $A_1$
and $A_2$. By corollary \ref{cor:matrixCTFE}, $g$ must be
constant. So
\begin{equation}\label{eq:solution}
f(x) = c \cdot x + b
\end{equation}
for some $c \in \mathbb{R}^n$ and $b \in \mathbb{R}$. Direct
substitution in (\ref{eq:cauchyveceq}) shows that
(\ref{eq:solution}) is a solution if and only if $b=0$.
\end{proof}

\begin{expl}
\emph{Let $\mathbf{K} = \{(x,y) \in \mathbb{R}^2 : |x| + |y| \leq
1\} $, ($K$ is the unit ball in $\ell_1 (\mathbb{R}^2)$), and let}
\begin{displaymath}
a_1(x,y) = \left(\frac{1}{2}x + \frac{1}{4}\sin y, \frac{1}{3}y
\right)
\end{displaymath}
\begin{displaymath}
a_2(x,y) = \left(\frac{1}{2}x - \frac{1}{4}\sin y, \frac{2}{3}y
\right)
\end{displaymath}
\emph{Denote by $\|\cdot\|_{\ell_1}$ the norm in $\ell_1
(\mathbb{R}^2)$. A straightforward computation shows that for} $i
= 1,2$
\begin{displaymath}
\|a_i (x,y)\|_{\ell_1} \leq \frac{11}{12}\|(x,y)\|_{\ell_1}
\end{displaymath}
\emph{and this shows that the $a_i$'s are maps in $K$ with an
attractor 0. The differentials of these maps are given by}
\begin{displaymath}
Da_1(x,y)= \left(%
\begin{array}{cc}
  \frac{1}{2} & \frac{1}{4} \cos y \\
  0 & \frac{1}{3} \\
\end{array}%
\right)
\end{displaymath}
\begin{displaymath}
Da_1(x,y)= \left(%
\begin{array}{cc}
  \frac{1}{2} & -\frac{1}{4} \cos y \\
  0 & \frac{2}{3} \\
\end{array}%
\right)
\end{displaymath}
\emph{Having all of the conditions of theorem
(\ref{thm:cauchyveceq}), we assert that the Cauchy type functional
equation}
\begin{displaymath}
f(x, y) - f\left(\frac{1}{2}x + \frac{1}{4}\sin y, \frac{1}{3}y
\right) - f\left(\frac{1}{2}x - \frac{1}{4}\sin y, \frac{2}{3}y
\right) = 0
\end{displaymath}
\emph{has only}
\begin{displaymath}
 f(x) = c \cdot x .
\end{displaymath}
\emph{as $C^1$ solutions.}
\end{expl}
\begin{expl}
\emph{Let} $\mathbf{K} = \{(x_1,x_2)^T \in \mathbb{R}^2 : {1 \over
2} \leq x_1^2 + x_2^2 \leq 1\}$, $\alpha = {\pi \over 3}$ ,
$$L_{\alpha} = \left(%
\begin{array}{cc}
  \cos (\alpha) & -\sin (\alpha) \\
  \sin (\alpha) & \cos (\alpha) \\
\end{array}%
\right)$$ \emph{and $R_{\alpha} = L_{\alpha}^T$. Let $\theta$
denote the angle between the positive $x$ axis and the line that
connects the point $(x_1,x_2)^T$ to the origin. Let $C_1(r,
\theta), C_2(r, \theta)$ be smooth functions that are periodic
with period ${\pi \over 3}$ in the second variable. Then}
$$f(x_1, x_2) = C_1(x_1^2 + x_2^2, \theta)x_1 + C_2(x_1^2 + x_2^2, \theta)x_2$$
\emph{is a solution to}
$$f(x) - f(L_{\alpha} x) - f(R_{\alpha} x)  = 0\verb"  ,  " x \in \mathbf{K} .$$
\end{expl}
In the above example, two conditions from theorem
(\ref{thm:cauchyveceq}) were violated: the eigenvalues are not
positive and the dynamical system generated by $L_\alpha$ and
$R_\alpha$ has no weak attractor. It would be interesting to find
a connection between the condition on the eigenvalues of the
differentials and the existence of a weak attractor. Perhaps
theorem (\ref{thm:cauchyveceq}) can be refined in such a way that
only conditions on the eigenvalues are given.

Theorem \ref{thm:cauchyveceq} was proved for general equations in
which appear general maps $a_i$, at the price of being able to
deal with compact domains only. But for a restricted family of
maps $a_i$ we can actually prove the ``uniqueness" of solutions to
the Cauchy type functional equation in the entire space
$\mathbb{R}^n$.

\begin{thm}\label{thm:affineCTFE}
Let $A_1$ and $A_2$ be two commuting, positive definite
(symmetric) $n \times n$ matrices and let $b_1, b_2 \in
\mathbb{R}^n$. Define for any $x \in \mathbb{R}^n$
\begin{displaymath}
T_i x = A_i x + b_i \hspace{5 mm} , \hspace{5 mm} i = 1,2
\end{displaymath}
All $C^1$ solutions $f:\mathbb{R}^n \rightarrow \mathbb{R}$ of the
Cauchy type functional equation
\begin{equation}\label{eq:entirespace}
f(T_1 x + T_2 x) = f(T_1 x) + f(T_2 x) \hspace{5 mm} , \hspace{5
mm} x \in \mathbb{R}^n
\end{equation}
are of the form
\begin{displaymath}
 f(x) = c \cdot x
\end{displaymath}
for some constant vector $c \in \mathbb{R}^n$.
\end{thm}
\begin{proof}
Let $f \in C^1$ be a solution of \ref{eq:entirespace}. Define a
new variable
\begin{displaymath}
y = S x \equiv (A_1 + A_2)x + b_1 + b_2
\end{displaymath}
then we may rewrite equation \ref{eq:entirespace} as
\begin{equation}\label{eq:entirespace'}
f(y) = f(T_1 S^{-1} y) + f(T_2 S^{-1} y) \hspace{5 mm} , \hspace{5
mm} y \in \mathbb{R}^n
\end{equation}
Note that $T_1 S^{-1} y = A_1 ((A_1 + A_2)^{-1} (y - b_1 -b_2) ) +
b_1$, so we introduce a matrix $B_1$
\begin{displaymath}
B_1 = A_1 (A_1 + A_2)^{-1}
\end{displaymath}
and a vector $d_1 \in \mathbb{R}^n$
\begin{displaymath}
d_1 = B_1 (-b_1 -b_2) + b_1
\end{displaymath}
to obtain the convenient form $T_1 S^{-1} y = B_1 y + d_1$.
Similarly, $T_2 S^{-1} y = B_2 y + d_2$, and we re-write
(\ref{eq:entirespace'}) as
\begin{equation}\label{eq:entirespace''}
f(y) = f(B_1 y + d_1) + f(B_2 y + d_2) \hspace{5 mm} , \hspace{5
mm} y \in \mathbb{R}^n
\end{equation}
>From the definitions it follows that $B_1 + B_2 = I$, and that all
the eigenvalues of $B_1, B_2$ are strictly between $0$ and $1$.
Being symmetric, the $B_i$'s are diagonalizable, thus there exists
a $\gamma < 1$ such that $\|B_i y\|_{\ell_2} \leq \gamma
\|y\|_{\ell_2}$ for $i = 1,2$.

Introduce the notation $\delta_i (y) = B_i y + d_i$ and
$\tilde{d_i} = \sum_{k = 0}^\infty B_i^k d_i$ , $i = 1,2$. For any
$y_1, y_2 \in \mathbb{R}^n$ we have $\|\delta_i (y_1) - \delta_i
(y_2)\|_{\ell_2} \leq \gamma\|y_1 - y_2\|_{\ell_2}$. On the other
hand
\begin{displaymath}
\delta_i (\tilde{d_i})  = B_i (\sum_{k = 0}^\infty B_i^k d_i) +
d_i = \sum_{k = 0}^\infty B_i^k d_i = \tilde{d_i}
\end{displaymath}
This means that for any point in $z \in \mathbb{R}^n$ the orbit
\begin{displaymath}
(z, \delta_i (z), \delta_i^2 (z), \ldots )
\end{displaymath}
converges exponentially to $\tilde{d_i}$. Now let $N$ be a
positive integer such that
\begin{displaymath}
N > \frac{\|\tilde{d_1} - \tilde{d_2}\|_{\ell_2}}{1-\gamma} \quad
.
\end{displaymath}
For each $m \geq N$ define
\begin{displaymath}
K_m = \overline{B}(\tilde{d_1},m) \cup \overline{B}(\tilde{d_2},m)
\end{displaymath}
where $\overline{B}(x,r)$ denotes the closed ball centered at $x$
with radius $r$. We note that the condition on $N$ insures that if
$x \in \overline{B}(\tilde{d_i},m)$, $i = 1,2$, and $i\neq j =
1,2$ then
\begin{displaymath}
\|\delta_j (x) - \tilde{d_j}\| \leq \gamma \|x-\tilde{d_j}\| \leq
\gamma \left(\|x-\tilde{d_i}\|+\|\tilde{d_i}-\tilde{d_j}\|\right)
\leq m
\end{displaymath}
It turns out that $K_m$ is compact, connected and
$\delta$--invariant, so for each $m$ we apply theorem
\ref{thm:cauchyveceq} with $K = K_m$ to infer that
\begin{displaymath}
f(y) = c_m \cdot y
\end{displaymath}
for all $y \in K_m$. But $\{K_m\}$ is an increasing sequence of
sets whose union is $\mathbb{R}^n$, so there is some $c$ such that
$c_m = c$ for all $m \geq N$. This shows what we claimed above.
\end{proof}

\pagebreak

\section{Unique solvability}\label{sec:UniSovCTFE} In the preceding section we used
a maximum principle to assert that, under some appropriate
conditions, the only solutions of the homogeneous equation
\begin{equation}\label{eq:cauchytypefunceq}
f(x) - \sum_{i=1}^N a_i (x) f(\delta_i (x)) = 0
\end{equation}
are constants. This clearly implies that, under the same
conditions, if the following non-homogeneous equation
\begin{equation}\label{eq:NHcauchytypefunceq}
f(x) - \sum_{i=1}^N a_i (x) f(\delta_i (x)) = h(x)
\end{equation}
has two solutions $f_1$ and $f_2$, then $f_1 = f_2 + C$ for some
constant $C$. In this section we shall also concern ourselves with
the \emph{solvabilty}, as well the uniqueness of solutions, of
functional equations of the type (\ref{eq:NHcauchytypefunceq}).

\begin{thm}\label{thm:uniquesolvability}
Let $(X,\delta)$ be a compact dynamical system. For $i = 1, \ldots
, N$, let $a_i : X \rightarrow \mathbb{R}$ be non-negative,
continuous functions such that
\begin{equation}\label{eq:coeffcond}
\forall x \in X . \sum_{i=1}^N a_i (x) \leq 1 .
\end{equation}
Define the guiding sets
\begin{displaymath}
\Lambda_i = \{x \in X : a_i (x) = 0\} .
\end{displaymath}
Assume that there is in $X$ a $\Lambda$-weak attractor $x_0$, and
that \newline $x_0 \in \{x \in X : \sum_{i=1}^N a_i (x) < 1\}$.
Then for any $h \in C(X)$ the functional equation
(\ref{eq:NHcauchytypefunceq}) has a unique solution $f \in C(X)$.
\end{thm}
\begin{rem} \emph{ This theorem was essentially proved by Paneah in
\cite{Pa03NC}, (Theorem 3). There $X$ was the interval $I=[-1,1]$
and the existence of an attractive set in $\partial I$ was a
consequence of explicit assumptions on $\delta$. The proof we give
is a modification of the proof given in \cite{Pa03NC}.}
\end{rem}

\begin{proof}
Define a linear operator $A : C(X) \rightarrow C(X)$ by
\begin{displaymath}
A f = \sum_{i=1}^N a_i \cdot f \circ \delta_i
\end{displaymath}
It is enough to prove that \footnote{In this proof, $\|\cdot\|$
will denote both the sup norm on $C(X)$ and the operator norm on
$L\left(C(X)\right)$, the space of bounded linear operators on
$C(X)$.}
\begin{equation}\label{eq:normOfA^m}
\exists m \in \mathbb{N} . \| A^m \| < 1 .
\end{equation}
Indeed, if this is the case, then the operator
\begin{displaymath}
f \mapsto f - Af
\end{displaymath}
is invertible \footnote{See \cite{Pa04DAIG} for a concise proof of
this fact.}, and this is exactly the content of the theorem. We
shall prove \ref{eq:normOfA^m} by a series of lemmas.
\begin{lem}
Let $T : C(X) \rightarrow C(X)$ be a positive linear operator.
Then $\|T\| = \|T {\bf 1} \|$.
\end{lem}
\begin{proof}
Let $f \in C(X)$ be of norm $1$. Then ${\bf 1} - f \geq 0$ thus
$T({\bf 1} - f) \geq 0$ or $T {\bf 1} \geq Tf$. Similarly, $-T
{\bf 1} \leq Tf$. This clearly implies $\|T {\bf 1} \| \geq \|Tf
\|$, and the lemma follows.
\end{proof}
For every $n \in \mathbb{N}$, define a continuous function $g_n$
on $X$ by
\begin{displaymath}
g_n (x) = \left( A^n {\bf 1} \right) (x) .
\end{displaymath}
Note that $A$ is a positive operator. By the above lemma, it
suffices to show that
\begin{equation}\label{eq:sufficesgn}
\exists m \in \mathbb{N} . \| g_m \| < 1 .
\end{equation}
Let's take a closer look at the functions $g_n$.
\begin{lem}
Explicitly, for $n \geq 2$, $g_n$ is given by
\begin{equation}\label{eq:g_n}
g_n (x) = \sum_{i_1, \ldots, i_n} a_{i_n} (x) \cdot a_{i_{n-1}}
(\delta_{i_n} (x)) \cdots a_{i_1} (\delta_{i_2} \circ \cdots \circ
\delta_{i_n} (x))
\end{equation}
where the sum is over all multi--indices $(i_1, \ldots, i_n) \in
\{1, \ldots, N \}^n$.
\end{lem}
\begin{proof}
We use induction.
\begin{displaymath}
g_1 (x) = \sum_{i=1}^N a_i (x)
\end{displaymath}
and
\begin{displaymath}
g_2 (x) = \left(A g_1 \right) (x) = \sum_{j=1}^N a_j (x) \cdot
\sum_{i = 1}^N a_i (\delta_j (x))
\end{displaymath}
and this is (\ref{eq:g_n}) for $n=2$. Now let $n > 2$.
\begin{eqnarray*}
g_n (x) &=& \sum_{i_n = 1}^N a_{i_n} (x) g_{n-1} (\delta_{i_n} (x)) \\
        &=& \sum_{i_n = 1}^N a_{i_n} (x) \sum_{i_1, \ldots ,i_{n-1}}
            a_{i_{n-1}} (\delta_{i_n} (x)) \cdots a_{i_1} (\delta_{i_2} \circ \cdots \circ
            \delta_{i_n} (x)) \\
        &=& \sum_{i_1, \ldots, i_n} a_{i_n} (x) \cdot a_{i_{n-1}}
            (\delta_{i_n} (x)) \cdots a_{i_1} (\delta_{i_2} \circ \cdots \circ
            \delta_{i_n} (x))
\end{eqnarray*}
and (\ref{eq:g_n}) is proved.
\end{proof}
\begin{lem}\label{lem:monotone g_n}
For all $x \in X$, if $n<k$ then $g_k (x) \leq g_n (x)$.
\end{lem}
\begin{proof}
By the previous lemma,
\begin{eqnarray*}
g_n (x) &=& \sum_{i_2, \ldots , i_n} a_{i_n} (x) \cdots a_{i_2}
            (\delta_{i_3} \circ \cdots \circ \delta_{i_n} (x)) \sum_{i_1} a_{i_1} (y_{i_2, \ldots
            , i_n})\\
        &\leq& g_{n-1} (x)
\end{eqnarray*}
where we have denoted $y_{i_2, \ldots , i_n} = \delta_{i_2} \circ
\cdots \circ \delta_{i_n} (x)$ and used (\ref{eq:coeffcond}) for
the inequality.
\end{proof}

The following lemma will make the conclusion of the theorem quite
clear.
\begin{lem}
For any $x \in X$ there exists a positive integer $m(x)$ such that
\begin{equation}\label{eq:mx}
g_{m(x)} (x) < 1 .
\end{equation}
\end{lem}
\begin{proof}
Fix $x \in X$. Since $\{x \in X : \sum_{i=1}^N a_i (x) < 1\}$ is
open, there exists an open neighborhood $V$ of $x_0$ that is
contained in $\{x \in X : \sum_{i=1}^N a_i (x) < 1\}$. $x_0$ is a
$\Lambda$-weak attractor, so there exists a $\Lambda$-proper orbit
\begin{displaymath}
\left(x, \delta_{j_n} (x) , \delta_{j_{n-1}} (\delta_{j_n} (x)),
\ldots , \delta_{j_2} \circ \cdots \circ \delta_{j_n} (x) \right)
\end{displaymath}
emanating from $x$ and terminating in $V$. This means that
$\delta_{j_2} \circ \cdots \circ \delta_{j_n} (x) \in V$. Now, as
we have noted before,
\begin{displaymath}
g_n (x) = \sum_{i_2, \ldots , i_n} a_{i_n} (x) \cdots a_{i_2}
            (\delta_{i_3} \circ \cdots \circ \delta_{i_n} (x)) \sum_{i_1} a_{i_1} \left (\delta_{i_2} \circ
            \cdots \circ \delta_{i_n} (x) \right) .
\end{displaymath}
We may write the right hand side as
\begin{eqnarray*}
a_{j_n} (x) \cdots a_{j_2} (\delta_{j_3} \circ \cdots \circ
\delta_{j_n} (x)) \sum_{i_1} a_{i_1} \left(\delta_{j_2} \circ
\cdots \circ \delta_{j_n} (x) \right) + \\
 \sum_{i_2, \ldots , i_n \neq j_2, \ldots , j_n} a_{i_n} (x) \cdots a_{i_2}
            (\delta_{i_3} \circ \cdots \circ \delta_{i_n} (x)) \sum_{i_1} a_{i_1} (\delta_{i_2} \circ
\cdots \circ \delta_{i_n} (x))
\end{eqnarray*}
But
\begin{displaymath}
\sum_{i_2, \ldots , i_n} a_{i_n} (x) \cdots a_{i_2}
            \left(\delta_{i_3} \circ \cdots \circ \delta_{i_n} (x) \right) =
            g_{n-1} (x)
\end{displaymath}
and because $\left(x, \delta_{j_n} (x) , \delta_{j_{n-1}}
(\delta_{j_n} (x)), \ldots , \delta_{j_2} \circ \cdots \circ
\delta_{j_n} (x) \right)$ is $\Lambda$-proper we have that
$a_{j_n} (x) \cdots a_{i_2} (\delta_{j_3} \circ \cdots \circ
\delta_{j_n} (x)) \neq 0.$ Moreover, $\delta_{j_2} \circ \cdots
\circ \delta_{j_n} (x) \in V$, so $\sum_{i_1} a_{i_1}
\left(\delta_{j_2} \circ \cdots \circ \delta_{j_n} (x) \right) <
1$ and thus
\begin{displaymath}
g_n (x) < g_{n-1} (x) \leq 1.
\end{displaymath}
Taking $m(x) =n$ the proof is complete.
\end{proof}
We are now in a position to finish the proof of the theorem. For
every $x \in X$ there is an $m(x)$ such that
\begin{displaymath}
g_{m(x)} (x) < 1 .
\end{displaymath}
Since $g_{m(x)}$ is continuous, there is a neighborhood $V_x$ of
$x$ where
\begin{displaymath}
\forall y \in V_x . g_{m(x)} (y) < 1 .
\end{displaymath}
The neighborhoods $\{V_x \}_{x \in X}$ form an open covering of
the space $X$, and therefore, by compactness of $X$, there is a
finite sub-covering $\{V_{x_1}, \ldots , V_{x_k} \}$. Denote $m_j
= m(x_j)$ , $j=1, \ldots, k$, and put $m = \max \{m_1, \ldots ,
m_k \}$. Then for any $y \in X$ there is a $j \in \{1, \ldots , k
\} $ such that $y \in V_{x_j}$. So $g_{m_j} (y) < 1$. But by lemma
\ref{lem:monotone g_n}
\begin{displaymath}
g_m (y) \leq g_{m_j} (y) < 1
\end{displaymath}
so that the inequality  $g_m (y) <1$ holds for all $y \in X$.
Consequently
\begin{displaymath}
\| g_m \| < 1 \,\, ,
\end{displaymath}
and this completes the proof of theorem
\ref{thm:uniquesolvability}.
\end{proof}

\pagebreak

\section{The initial value problem for a $\mathcal{P}$-configuration}
In the previous sections we dealt with rather general dynamical
systems and functional equations. Now we will concentrate on a
very specific family of dynamical systems and their corresponding
Cauchy type functional equations. In fact, we shall prove a
necessary and sufficient condition for the existence of a unique
solution $f \in C^2(I)$ to the problem
\begin{eqnarray}\label{eq:Pconfeq1}
f(t) - f(\delta_1 (t)) - f(\delta_2 (t)) & = & h(t) \quad , \quad
t \in I \\ \label{eq:Pconfeq2} f'(c) & = & \mu
\end{eqnarray}
where $I = [a, b]$, $c \in (a,b)$, $\mu$ is some real number, $h
\in C^2$ satisfies $h(a) = h(b)$, and $\delta_1, \delta_2$ form a
$\mathcal{P}$-configuration in $I$. This problem is of great
importance for us for two reasons: 1) it is equivalent to a
boundary value problem which we treat in chapter 4, and 2)
``historically" the dynamical system in this problem is the origin
of the theory of guided dynamical systems. The history of this
problem can be found in Paneah's papers \cite{Pa97} -
\cite{Pa04DAIG}, where certain conditions for unique solvability
of the problem (\ref{eq:Pconfeq1})-(\ref{eq:Pconfeq2}) are proved.

\subsection{Definition of a $\mathcal{P}$-configuration}
Let $I = [a, b]$ be a fixed closed interval in $\mathbb{R}$, $c
\in (a,b)$, and let $\delta_1, \delta_2: I \rightarrow I$ be two
$C^2$ maps satisfying the following conditions:
\begin{equation}\label{eq:Pconf1}
\delta_1'(t) + \delta_2'(t) = 1 \quad , \quad t \in I \,\, ;
\end{equation}
\begin{equation}\label{eq:Pconf2}
\delta_i'(t) \geq 0 \quad , \quad t \in I, i = 1,2 \,\, ;
\end{equation}
\begin{equation}\label{eq:Pconf3}
\delta_2(a) = a, \quad \delta_2(b) = \delta_1(a) = c, \quad
\delta_1(b) = b \,\, .
\end{equation}
If all these assumptions hold, then the maps $\delta_1$ and
$\delta_2$ are said to form a $\mathcal{P}$-configuration in $I$.
We introduce the guiding sets
\begin{displaymath}
\Lambda_1 = \{t \in I | \delta_1'(t) = 0\}
\end{displaymath}
and
\begin{displaymath}
\Lambda_2 = \{t \in I | \delta_2'(t) = 0\} .
\end{displaymath}

\subsection{Generalized $\mathcal{P}$-configuration} At the same
cost of proving the necessary and sufficient conditions for unique
solvability of (\ref{eq:Pconfeq1})-(\ref{eq:Pconfeq2}), we may
prove the same type of theorem for a class of a equations that is
a little more general. To this end, we make the following
definitions.

\begin{figure}[t]
\centering
\includegraphics[scale=0.5]{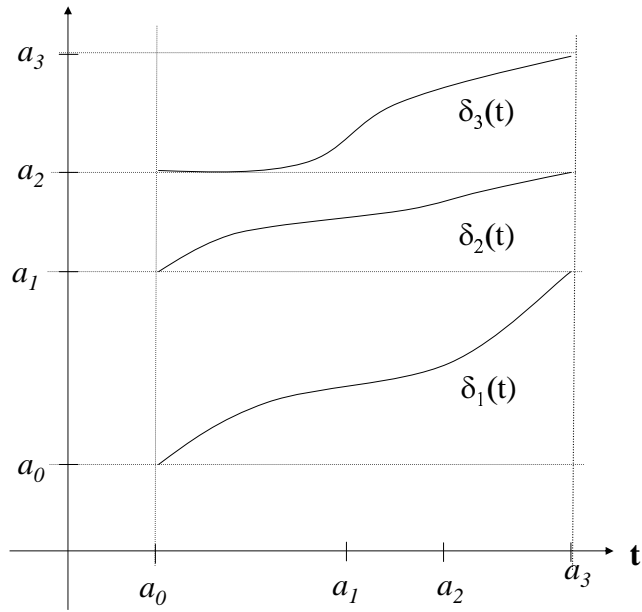}
\caption{Generalized $\mathcal{P}$-configuration.}
\label{fig:Pconf}
\end{figure}

Let $a_0 < a_1 < \ldots < a_N$ be $N+1$ points in $\mathbb{R}$.
Define $I = [a_0, a_N]$. Let $\delta_1, \ldots, \delta_N$ be $C^2$
functions such that $\delta_i$ maps $I$ onto $[a_{i-1}, a_i]$, for
$i = 1, \ldots, N$. Assume that
\begin{equation}\label{eq:GPconf1}
\sum_{i=1}^N \delta_i'(t) = 1 \quad , \quad t \in I \,\,,
\end{equation}
\begin{equation}\label{eq:GPconf2}
\delta_i'(t) \geq 0 \quad , \quad t \in I \,,\,\,
i=1,\ldots,N\,\,,
\end{equation}
\begin{equation}\label{eq:GPconf3}
\delta_i(a_0) = a_{i-1}, \, \,  \delta_i(a_N) = a_i \quad , \quad
i=1,\ldots,N\,\,.
\end{equation}
We say that the maps $\delta_1, \ldots, \delta_N$ generate a
\emph{generalized $\mathcal{P}$-configuration} in $I$. For $i = 1,
\ldots, N$, introduce the guiding sets
\begin{displaymath}
\Lambda_i = \{t \in I \, | \, \delta_i'(t) = 0 \} \, .
\end{displaymath}
See figure \ref{fig:Pconf}.

Our aim now will be to prove a necessary and sufficient condition
for the existence of a unique solution $f$ to the following
problem:
\begin{eqnarray}\label{eq:GPconfeq1a}
f(t) - \sum_{i=1}^N f(\delta_i (t)) & = & h(t) \quad , \quad t \in
I
\\ \label{eq:GPconfeq2b} f'(c) & = & \mu \quad ,
\end{eqnarray}
where the point $c \in I$ and the number $\mu$ are given, and $h$
is an arbitrary $C^2$ function satisfying $h(a_0) = h(a_N)$.

\subsection{Some preliminary results in functional analysis and
preparations} \label{subseq:perliminary}
In this section we shall
use without explanation results from functional analysis. Our
reference for facts regarding Fredholm operators and
Riesz-Schauder theory is \cite{FuncAnalysis}. Let us just recall
the following two facts:
\begin{enumerate}
  \item If $A$ is a Fredholm operator and $K$ is compact then $A +
  K$ is also Fredholm and $\textrm{ind}(A+K) = \textrm{ind}(A)$ \footnote{For a Fredholm operator
  $A$ we denote by $\textrm{ind}(A)$ the index of $A$}.
  \item If $A : \textbf{V} \rightarrow \textbf{V}'$ and $B : \textbf{V}' \rightarrow
  \textbf{V}''$ are Fredholm, then $BA$ is also Fredholm and $\textrm{ind}(BA) = \textrm{ind}(A) + \textrm{ind}(B)$
\end{enumerate}

Before proceeding it is worth noting that the idea to use
Riesz-Schauder theory in this problem is due to Paneah and was
introduced in the papers cited above. However, as is the case in
papers many times, this idea was not explained in great detail.
Therefore, the rest of this subsection is devoted to making the
necessary preparations that will justify the use we shall make
later on of Paneah's idea.

Fix some point $c \in [a_0,a_N]$. We introduce the function spaces
\footnote{These are Banach spaces when equipped with the usual
norm. For example, $\|f\|_{{}_{\textbf{X}}} = \sup_I|f| +
\sup_I|f'| + \sup_I|f''|$, etc.}
\begin{eqnarray*}
\textbf{X} &=& \left\{\varphi \in C^2(I) \, | \, \sum_{i=1}^{N-1}\varphi(a_i) = \varphi'(c) = 0 \right\} \\
\textbf{Y} &=& \left\{\psi \in C^2(I) \, | \, \psi(a_0) =
\psi(a_N) = 0 \right\} \\ \textbf{W} &=& \left\{\xi \in C^1(I) \,
| \, \xi(c) = 0 \right\}
\end{eqnarray*}
and
\begin{displaymath}
\textbf{Z} = \left\{\omega \in C^1(I) \, | \, \int_{a_0}^{a_N}
\omega(t) dt = 0 \right\} .
\end{displaymath}
Define $B_0 \in L(\textbf{X},\textbf{Y})$, $B_1 \in
L(\textbf{W},\textbf{Z})$ and $B_2 \in L\left(C(I)\right)$ by
\begin{eqnarray*}
(B_0 f)(t) &=& f(t) - \sum_{i=1}^N f\left(\delta_i (t) \right) \quad \\
 (B_1 g)(t) &=& g(t) - \sum_{i=1}^N \delta_i'(t)
g\left(\delta_i(t) \right)
\end{eqnarray*}
and
\begin{displaymath}
(B_2 h)(t) = h(t) - \sum_{i=1}^N {\delta_i'}^2 h(\delta_i(t)) -
\sum_{i=1}^N \delta_i'' \int_{c}^{\delta_i(t)} h(s) ds \, .
\end{displaymath}
An easy check shows that these operators are bounded (with respect
to the standard norms of these spaces) and that they map into the
right spaces. For example, if $f \in \textbf{X}$, then
\begin{eqnarray*}
(B_0 f)(a_0) &=& f(a_0) - \sum_{i=1}^N f\left(\delta_i
(a_0)\right)
\\ &=& f(a_0) - f(a_0) - \sum_{i=1}^{N-1} f(a_i) = 0
\end{eqnarray*}
so $(B_0 f) (a_0) = 0$ and $(B_0 f) (a_N) = 0$ is shown in a
similar manner, thus $(B_0 f) \in \textbf{Y}$. There are four
different invertible bounded linear operators $\textbf{X}
\rightarrow \textbf{W}$, $\textbf{Y} \rightarrow \textbf{Z}$,
$\textbf{W} \rightarrow C(I)$ and $\textbf{Z} \rightarrow  C(I)$
representing differentiation. Let us make a convenient abuse of
notation by denoting all of these operators by $D$.

Differentiating equation (\ref{eq:GPconfeq1a}) once and twice
gives
\begin{eqnarray}\label{eq:PconfeqD}
\left((B_0 f)(t)\right)' &=& (B_1 f')(t) \\ \nonumber &=& f'(t) -
\sum_{i=1}^N \delta_i'(t)\cdot f'\circ\delta_i(t) \\ \nonumber
 &=& h'(t)
\end{eqnarray}
and (if $f \in \textbf{X}$)
\begin{eqnarray}\label{eq:PconfeqD2}
\left((B_0 f)(t)\right)'' &=& (B_2 f'')(t) \\ \nonumber &=& f''(t)
- \sum_{i=1}^N {\delta_i'}^2 f''(\delta_i(t)) - \sum_{i=1}^N
\delta_i'' f'(\delta_i(t)) \\ \nonumber  & = & f''(t) -
\sum_{i=1}^N {\delta_i'}^2 f''(\delta_i(t)) - \sum_{i=1}^N
\delta_i'' \int_{c}^{\delta_i(t)} f''(s) ds
\\ \nonumber
 &=& h''(t) \, .
\end{eqnarray}
(We used the fact that $f'(c) = 0$). From this it follows that
\begin{equation} \label{eq:DBBD1}
D B_0 = B_1 D
\end{equation}
and
\begin{equation} \label{eq:DBBD2}
D B_1 = B_2 D \, .
\end{equation}
\begin{lem}\label{lem:equiv}
If one of $B_0, B_1$ or $B_2$ is injective (surjective), then all
of $B_0, B_1$ and $B_2$ are injective (surjective). If one of
$B_0, B_1$ or $B_2$ is Fredholm, then all of $B_0, B_1$ and $B_2$
are Fredholm and  \, $\emph{\textrm{ind}}B_0 =
\emph{\textrm{ind}}B_1 = \emph{\textrm{ind}}B_2$.
\end{lem}
\begin{proof}
Assume, for instance, that $\textrm{Ker}B_0 = \{0\}$. Then
$\textrm{Ker}D B_0 = \{0\}$, whereas $\textrm{Ker}B_1 D = D^{-1}
(\textrm{Ker}B_1)$. From (\ref{eq:DBBD1}) we infer that
$\textrm{Ker}B_1 = \{0\}$. The rest of the first statement is
proved in a similar manner.

Abusing our notation a little more we may write, e.g.,
\begin{displaymath}
B_0 = D^{-1}B_1 D
\end{displaymath}
Taking into account the second fact that we cited above
\footnote{and also the fact that $D$ is bounded and invertible on
the relevant spaces} this shows that $B_0$ is Fredholm if and only
if $B_1$ is, and that their indices agree, since
\begin{displaymath}
{\rm ind}(B_0) = {\rm ind}\left(D^{-1}\right) + {\rm ind}(B_1) +
{\rm ind}(D) = 0 + {\rm ind}(B_1) + 0 = {\rm ind}(B_1)
\end{displaymath}
\end{proof}

\subsection{The initial value problem}\label{subsec:initvalprob}
\begin{thm}\label{thm:Pconf1}
$\emph{\textrm{Im}}B_0 = \textbf{Y}$ if and only if $(I, \delta,
\Lambda)$ is $\Lambda$-minimal. When this is the case, $B_0$ is an
isomorphism.
\end{thm}
\begin{proof}
Let us begin by showing necessity. Assume that $(I, \delta,
\Lambda)$ is not $\Lambda$-minimal. We have to show that
$\textrm{Im}B_0 \neq \textbf{Y}$. By lemma \ref{lem:equiv}, it is
enough to show that $B_1$ is not surjective.

By proposition \ref{prop:miniffno}, there exists a closed,
non-empty $(\Lambda, \delta)$-invariant set $A \subsetneq I$. Let
$G$ be a $C^1$ function such that
\begin{displaymath}
\int_{a_0}^{a_N} G(t) dt = 0
\end{displaymath}
and $G \big| _A \equiv 1$. Attempting to arrive at a contradiction
we assume that $F \in \textbf{W}$ is a solution to the equation
\begin{equation} \label{eq:BFG}
B_1 F = G \, .
\end{equation}
Denote $M = \max_{t \in I} |F(t)|$. Define a linear operator $T:
C(I) \rightarrow C(I)$ by
$$T F = \sum_{i=1}^N \delta_i'\cdot F\circ\delta_i \quad .$$
By (\ref{eq:GPconf1}) and (\ref{eq:GPconf2}), $\|T\| \leq 1$. As
$A$ is $(\Lambda, \delta)$-invariant, we also have that for all
$k$, $(T^k G) \big|_A \equiv 1$. Fix some $t_0 \in A$, and let
$\textbf{I}$ denote the identity operator on $C(I)$.  Operating on
both sides of (\ref{eq:BFG}) with the operator $\textbf{I} + T +
T^2 + \ldots + T^n$ at the point $t_0$, and noting that$B_1 =
\textbf{I} - T$, we obtain
\begin{displaymath}
(\textbf{I} - T^{n+1})F(t_0) = (\textbf{I}+T+\dots+T^n)G(t_0)
\end{displaymath}
thus for all $n$ we have that
\begin{displaymath}
2M \geq |(\textbf{I}-T^{n+1})F(t_0)| =
|(\textbf{I}+T+\dots+T^n)G(t_0)| = n + 1
\end{displaymath}
a contradiction.

Following Paneah, the sufficiency will be established by proving
that:
\begin{enumerate}
  \item $\textrm{Ker}B_1 = \{0\}$
  \item $B_2$ is a Fredholm operator and $\textrm{ind}B_2 = 0$.
\end{enumerate}
Recall that lemma \ref{lem:equiv} translates these facts to the
invertibility of $B_0$.

\emph{Proof of 1}. Let $F \in \textbf{W}$ satisfy $B_1 F = 0$.
Note that $B_1 F = 0$ is precisely the functional equation studied
in section \ref{sec:Maxp}. The conditions on the maps in a
$\mathcal{P}$-configuration, and the existence of $\Lambda$-weak
attractor, (which is a trivial consequence of
$\Lambda$-minimality), all add up to the fact that $F$ and $(I,
\delta, \Lambda)$ satisfy the conditions of lemma
\ref{lem:lemmaxp}, and thus $F = const.$ But, being in
$\textbf{W}$, $F (c) = 0$, thus $F = 0$. This proves 1.

\emph{Proof of 2}. Define the operators $L,K: C(I) \rightarrow
C(I)$
\begin{displaymath}
(L F)(t) = \sum_{i=1}^N {\delta_i'}^2 F(\delta_i(t))
\end{displaymath}
and
\begin{displaymath}
(K F)(t) = \sum_{i=1}^N \delta_i'' \int_{c}^{\delta_i(t)} F(s) ds
\, \, .
\end{displaymath}
With this new notation we can decompose $B_2$ as $B_2 = \textbf{I}
- L - K$. Now, $\delta_1, \ldots , \delta_N \in C^2$, so the set
where at least two of the $\delta_i'$ are positive is non-empty.
But this set is exactly
\begin{displaymath}
\left\{t \in I \, | \, \sum_{i=1}^N {\delta_i'}^2 (t) < 1 \right\}
\end{displaymath}
and by the assumed  $\Lambda$-minimality this set contains a
$\Lambda$-weak attractor. We can now employ theorem
\ref{thm:uniquesolvability} to conclude that $\textbf{I} - L$ is
an invertible operator\footnote{In this work, a function (or
operator) is called \emph{invertible} if it is both injective and
surjective.}. 2 now follows from the fact that $K$ is a compact
operator, and from the first fact from functional analysis cited
at the beginning of \ref{subseq:perliminary}.
\end{proof}
\begin{rem}\emph{
For applications in partial differential equations it is worth
noting that the operator $B_0^{-1}$ is bounded if the operator
$B_0$ is invertible. This, of course, follows from Banach's open
mapping theorem.}
\end{rem}
Note that in the above proof for sufficiency we used the
$\Lambda$-minimality only to infer the existence of a
$\Lambda$-weak attractor in $\mathcal{A}=\{t \in I \, | \, \sum
{\delta_i'}^2 (t) < 1 \}$. This set $\mathcal{A}$ contains $\{t \,
| \, \forall i \, . \, \delta_i' (t) >0 \} = I \setminus \Lambda$.
Thus the existence of a $\Lambda$-weak attractor in $I \setminus
\Lambda$ is a sufficient condition for the solvability of the
equation $B_0 f = h$. But we have just shown that the solvability
of this problem implies that $(I, \delta, \Lambda)$ is
$\Lambda$-minimal! Thus we arrive at the very unexpected result:
\begin{prop}\label{prop:miniffwa}
In a $\mathcal{P}$-configuration $(I, \delta, \Lambda)$ the
following are equivalent:
\begin{enumerate}
  \item $(I, \delta, \Lambda)$ is $\Lambda$-minimal.
  \item There exists a $\Lambda$-weak attractor in $I \setminus \Lambda$.
\end{enumerate}
\end{prop}

Now we return to the problem (\ref{eq:GPconfeq1a}) -
(\ref{eq:GPconfeq2b}).

\begin{thm}\label{thm:initvalprob}
Let $(I, \delta, \Lambda)$ be a generalized
$\mathcal{P}$-configuration that has a $\Lambda$-weak attractor in
$I \setminus \Lambda$. Then for any $h \in C^2(I)$ with $h(a_0) =
h(a_N)$, and for any $\mu \in \mathbb{R}$, $c \in [a_0,a_N]$,
there exists a unique solution $f \in C^2(I)$ of the problem
\begin{eqnarray}\label{eq:GPconfeq1}
f(t) - \sum_{i=1}^N f(\delta_i (t)) & = & h(t) \quad , \quad t \in
I
\\ \label{eq:GPconfeq2} f'(c) & = & \mu
\end{eqnarray}
\end{thm}
\begin{rem}\emph{
Substituting $t = a_0$ and $t = a_N$ in (\ref{eq:GPconfeq1}) we
see, using the properties of the $\delta$'s, that if $f$ is a
solution to (\ref{eq:GPconfeq1}) then
$$\sum_{i=1}^{N-1}f(a_i) = h(a_0) = h(a_N) \, .$$}
\end{rem}

%

\begin{proof}
Let $h \in C^2(I)$ satisfy $h(a_0) = h(a_N)$. Define
\begin{displaymath}
\tilde{h}(t) = h(t) - h(a_0) \, .
\end{displaymath}
Using the notation introduced in \ref{subseq:perliminary}, we have
that $\tilde{h} \in \textbf{Y}$. By theorem \ref{thm:Pconf1},
there exists an $\tilde{f} \in \textbf{X}$ such that
\begin{displaymath}
\tilde{f}(t) - \sum_{i=1}^N \tilde{f}(\delta_i(t)) = \tilde{h}(t)
\quad , \quad t \in I \, .
\end{displaymath}
Put
\begin{displaymath}
f(t) = \tilde{f}(t) - \frac{h(a_0) + \mu C}{N-1} + \mu t
\end{displaymath}
where $C$ satisfies $\sum_{i=1}^{N}\delta_i (t) = t + C$.


Now $f'(c) = \mu$, and
\begin{eqnarray*}
(B_0 f)(t) &=& (B_0 \tilde{f})(t) - B_0 \left( \frac{h(a_0) + \mu C}{N-1} \right) + B_0 (\mu t) \\
\nonumber  &=& \tilde{h}(t) + h(a_0) + \mu C - \mu C\\
           &=& h(t) \, .
\end{eqnarray*}
Uniqueness follows from \ref{thm:Pconf1}.
\end{proof}
%

\chapter{Overdeterminedness of functional equations}\label{ch:OD}
The branch in mathematics that is concerned with functional
equations splits into two main sub-branches, dealing with two main
sub-classes of equations, namely ``functional equations in a
single variable''and ``functional equations in several
variables''. Up to now we have only considered equations that
belong to the first class. In this section we will address some
problems that lie on the borderline between these two classes.

Recall the classical Cauchy functional equation:
\begin{equation}\label{eq:ClassicalCauchyEq}
f(x + y) = f(x) + f(y)
\end{equation}
This is a functional equation in 2 variables. To \emph{solve} the
functional equation usually means : given a set $A \subseteq
\mathbb{R}^2$ and a class of functions $\mathcal{A}$, to find the
family of functions $\mathcal{F} \subseteq \mathcal{A}$ which
consists of all $f$ such that $f(x+y) = f(x) + f(y)$ for all
$(x,y) \in A$. Following Kuczma (\cite{KuczmaRestrictedDomains})
let us call $A$ the \emph{domain of validity}. For example, when
Cauchy first treated (\ref{eq:ClassicalCauchyEq}), he took
$\mathcal{A} = C(\mathbb{R})$, and showed that if the domain of
validity is taken to be $\mathbb{R}^2$ then the only solutions to
(\ref{eq:ClassicalCauchyEq}) are of the form $f(z) = \lambda z$.
It has been shown in various works (\cite{Sablik1},
\cite{Sablik2}, \cite{Aczel2}, \cite{Pa04OD} and the references
therein) that when some additional smoothness assumptions are
imposed on $f$ then even if the domain of validity is quite small
- the graph of an appropriate function, for example - the set of
solutions doesn't grow. Thus, using the terminology of Paneah
(\cite{Pa04OD}), we may say that the equation
\begin{equation*}
f(x + y) = f(x) + f(y) \hspace{5 mm} , \hspace{5 mm} (x,y) \in
\mathbb{R}^2
\end{equation*}
is \emph{overdetermined} (for the class of functions satisfying
these additional smoothness assumptions). For an explicit example,
consider the equation
\begin{equation*}
f(t) - f\left(\frac{t+1}{2}\right) - f\left(\frac{t-1}{2}\right) =
0 \hspace{5 mm} , \hspace{5 mm} t \in [-1,1]
\end{equation*}
This is nothing but the classical Cauchy equation with domain of
validity $\Gamma = \{\left( (t+1)/2, (t-1)/2 \right) \mid t \in
[-1,1]\}$. By theorem \ref{thm:cauchyveceq}, the only $C^1$
functions satisfying this equation are $f(z) = \lambda z$. Thus
theorem \ref{thm:cauchyveceq} may be interpreted as the assertion
that the equation
\begin{equation*}
f(x + y) = f(x) + f(y) \hspace{2 mm} , \hspace{3 mm} |x+y| \leq 1,
|x|,|y| \leq 1
\end{equation*}
is \emph{overdetermined} for functions for the class $C^1$. In
fact, note that in subsection \ref{subsec:CauchyRn} we proved that
the Cauchy equation in $\mathbb{R}^n$
\begin{displaymath}
f(x_1+y_1, \ldots, x_n+y_n) = f(x_1, \ldots, x_n) + f(y_1, \ldots,
y_n)
\end{displaymath}
is overdetermined for the class $C^1(\mathbb{R}^n,\mathbb{R})$

One is led to the following questions: (a) \emph{given a class of
functions $\mathcal{A}$, what is the ``smallest'' domain of
validity for which the solutions to (\ref{eq:ClassicalCauchyEq})
are only $f(z) = \lambda z$}, and : (b) \emph{given a domain of
validity, for what $\mathcal{A}$ does the set of solutions to
(\ref{eq:ClassicalCauchyEq}) remain $f(z) = \lambda z$?}

The above questions may be asked with regards to any functional
equation, and it is interesting in general to study how, given a
functional equation, the set of solutions changes when the domain
of validity and the class of functions considered are changed.
This direction of research attracted relatively little attention
during the years, and most of the efforts were put into Cauchy's
equation. Before we can continue, it is important to note that the
terminology we use is not standard. There is no way to escape
this, as practically every researcher in this field used different
terminology. M. Kuczma used the term \emph{functional equations on
restricted domains} to describe the general problem
(\cite{KuczmaRestrictedDomains}), while Acz\'{e}l and Dhombres
prefer \emph{conditional functional equations} (\cite{Aczel2}).
Synonyms for \emph{overdeterminedness} are
\emph{redundancy}(Introduced by Dhombres and Ger in papers cited
in \cite{KuczmaRestrictedDomains}) and in some places
\emph{addundancy}.

For most classical functional equations in 2 variables, the domain
of validity is usually taken to be some large, open set in
$\mathbb{R}^2$. In \cite{Pa04OD} Paneah proved for a sample of
classical functional equations that, under some smoothness
assumptions, their solution is already determined by the
functional equation holding on a much smaller domain of validity,
e.g., a one-dimensional sub-manifold in $\mathbb{R}^2$, and such
equations were called \emph{overdetermined}. In this chapter we
prove two results in this spirit.

\section{Overdeterminedness of Cauchy's functional equation} In
this subsection we shall show the overdeterminedness of the Cauchy
functional equation for \emph{continuous} functions. It must be
noted that this fact follows immediately from the results of M.
Lackovich, who showed in \cite{Lackovich} that if $f: [0,\infty)
\rightarrow \mathbb{R}$ is \emph{measurable} and satisfies
Cauchy's equation on the line $\{(at, bt) | t \in \mathbb{R} \}$
where $\log_a b \notin \mathbb{Q}$, then $f(z) = \lambda z$.
\footnote{Lackovich' result was in fact much more general, we are
only stating the consequence that is directly connected to our
work.} We note that if $\log_a b \in \mathbb{Q}$, then for any
continuous function $A(z)$ with a period $1$, the function
\begin{displaymath}
f(z) = z \cdot A(log_c (z))
\end{displaymath}
is a continuous solution to the Cauchy fuunctional equation on
that line, where $c$ is some number that satisfies $c^k = a$ and
$c^m = b$ for some integer $k$ and $m$.

Define
\begin{equation*}
\Gamma = \{(x,y) \in \mathbb{R}^2 : |x| + |y| = 1 \}
\end{equation*}
and
\begin{equation*}
\Gamma^* = \Gamma \setminus \{(x,x+1) : x \in [-1,0] \}
\end{equation*}

\begin{figure}[t]
\centering
\includegraphics[scale=0.5]{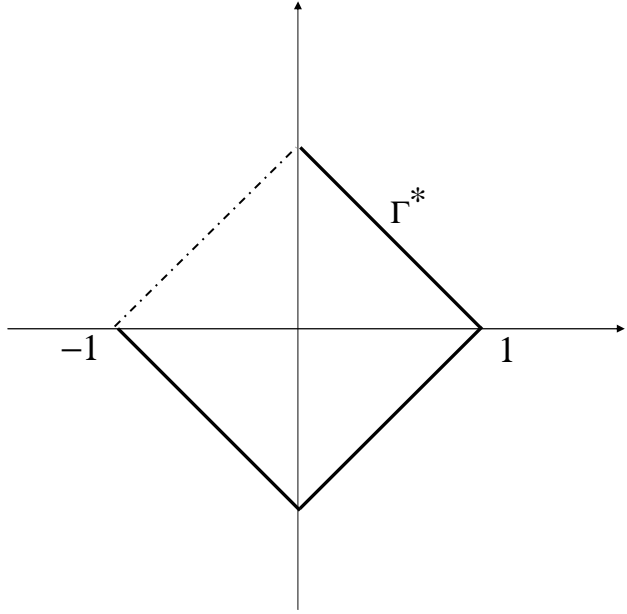}
\caption{The solutions of the Cauchy equation are determined on
$\Gamma^*$.}
\label{fig:cauchy}
\end{figure}

See figure \ref{fig:cauchy}. As mentioned above, if a continuous
function $f : [-1,1] \rightarrow \mathbb{R}$ satisfies the Cauchy
functional equation in the set $\{(x,y) : |x| + |y| \leq 1\}$ then
$f(z) = \lambda z$. We shall now show that the Cauchy functional
equation on the boundary of this set already determines the same
set of solutions. This result was first published by my
students\footnote{Ardazi, Kharash, Mamane and Zoabi.} in
\cite{SciTech}.
\begin{thm}\label{thm:ODofCE}
Let $f : [-1,1] \rightarrow \mathbb{R}$ be a continuous function
that satisfies the relation
\begin{equation*}
f(x+y) = f(x) + f(y) \hspace{2 mm} , \hspace{3 mm} (x,y) \in
\Gamma^*
\end{equation*}
Then $f(z) = \lambda z$ with a constant $\lambda$.
\end{thm}
Thus, the Cauchy equation in the square $\{(x,y) : |x| + |y| \leq
1\}$ is overdetermined for \emph{continuous} functions.
\begin{proof}
Fix the notation $I = [-1,1]$. Choosing parameterization for each
side of $\Gamma^*$ we arrive, after some simple manipulations (see
\cite{SciTech}) at the system of functional equations:
\begin{eqnarray}\label{eq:CauchySys}
  f\left( \frac{t-1}{2} \right) = \frac{f(t)-f(1)}{2} \hspace{2 mm} , \hspace{3
  mm}t \in I
  \\ \label{eq:CauchSys2}
  f\left( \frac{t+1}{2} \right) = \frac{f(t)+f(1)}{2} \hspace{2 mm} , \hspace{3
  mm}t \in I
\end{eqnarray}
Now introduce two maps $\alpha, \beta : I \rightarrow I$
\begin{displaymath}
\alpha(t) = \frac{t+1}{2}
\end{displaymath}
\begin{displaymath}
\beta(t) = \frac{t-1}{2}
\end{displaymath}
Clearly, $\alpha$ and $\beta$ are Lipschitz with constant
$\frac{1}{2}$, so by proposition \ref{prop:contracting} they
generate in $I$ the minimal dynamical system $(I,\{\alpha,
\beta\})$. In particular, the orbit-set of the point $1$ is dense
in $I$. To complete the proof, let us show that for every point $z
\in OS (1)$
\begin{displaymath}
f(z) = f(1) \cdot z.
\end{displaymath}
For $z = 1$ this is evident. Assume that $z_0 \in OS (1)$, and
that $f(z_0) = f(1) \cdot z_0$. Then
\begin{displaymath}
f(\alpha(z_0)) = \frac{f(z_0)+f(1)}{2}
\end{displaymath}
by \ref{eq:CauchySys}. But $f(z_0) = f(1) \cdot z_0$ so
\begin{displaymath}
\frac{f(z_0)+f(1)}{2} = \frac{f(1) \cdot z_0 + f(1)}{2} =
\frac{f(1) \cdot (z_0 + 1)}{2} = f(1) \cdot \alpha(z_0) \,\,.
\end{displaymath}
This means that $f(\alpha(z_0)) = f(1) \cdot \alpha(z_0)$.
Similarly, $f(\beta(z_0)) = f(1) \cdot \beta(z_0)$ and the theorem
follows.
\end{proof}

\section{A uniqueness/overdeterminedness theorem} In 1964
Acz\'{e}l \footnote{See \cite{AczelBooklet}.} proved the following
uniqueness theorem for a rather wide class of functional
equations:

\begin{thm}\label{thm:aczel}
Let $f_1, f_2 : I \rightarrow \mathbb{R} $ be continuous solutions
of the equation

\begin{equation}\label{eq:feaczel}
f(F(x,y)) = H[f(x),f(y),x,y] \verb"  ,  " (x,y) \in I^2
\end{equation}
where $I$ is an (open, closed, half-open, finite or infinite)
interval. Suppose that $F:I^2 \rightarrow I$ is continuous and
internal that is,

$$ min(x,y) < F(x,y) < max(x,y) \, \, \emph{\textrm{ if }} \, \,  x \neq y $$
and that either $u \mapsto H(u,v,x,y)$ or $v \mapsto H(u,v,x,y)$
are injections. Further, let $a,b \in I$ and

\begin{displaymath}
f_1(a) = f_2(a) \, \, \emph{\textrm{ and }} \, \, f_1(b) = f_2(b)
\, .
\end{displaymath}
Then
\begin{displaymath}
\forall x \in I . f_1(x) = f_2(x) \, \, .
\end{displaymath}
\end{thm}

This theorem motivated much work on uniqueness theorems and has
been improved several times. Theorems in the same spirit were
proved for different classes of $F$ and $H$ and for more general
spaces ($\mathbb{R}^2$,$\mathbb{R}^n$, topological vector spaces,
\ldots \footnote{\cite{Aczel2} contains references to these
developments.}). In this section we will prove a refinement of the
above theorem which serves at once both as a uniqueness theorem
for (\ref{eq:feaczel}) and as a proof that all of the equations
that belong to the class treated below are overdetermined.

\begin{figure}[t]
\centering
\includegraphics[scale=0.5]{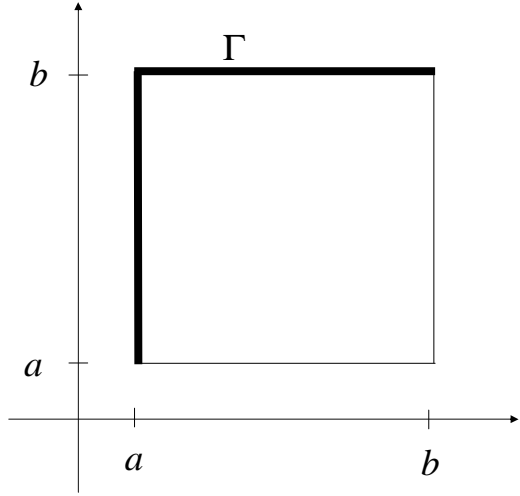}
\caption{The solution of the functional equation \ref{eq:feaczel} is
determined on $\Gamma$.}
\label{fig:OD}
\end{figure}

\begin{thm}\label{thm:uniqueOD}
Let $I = [a,b]$, $H : \mathbb{R} \times \mathbb{R} \times I \times
I \rightarrow \mathbb{R}$ any function and $F : I^2 \rightarrow I$
a continuous function that satisfies
\begin{itemize}
  \item $\forall x \neq y . |F(x,b)-F(y,b)|,|F(a,x)-F(a,y)|<|x-y|$
  \item $\exists x_0, y_0 . F(a,x_0) = a \quad {\rm and} \quad F(y_0,b) = b$
\end{itemize}
For any real $A$ and $B$ there exists at most one solution $f$ to
(\ref{eq:feaczel}) that satisfies the boundary conditions

\begin{equation}\label{eq:boundcond}
f(a) = A \quad f(b) = B .
\end{equation}
Moreover, if a function $f$ is a solution to (\ref{eq:feaczel})
satisfying (\ref{eq:boundcond}), then it is already determined by
the functional equation

\begin{equation}\label{eq:feongamma}
f(F(x,y)) = H[f(x),f(y),x,y] \quad , \quad (x,y) \in \Gamma
\end{equation}
 where $ \Gamma = ([a,b] \times \{b\} ) \cup (\{a\} \times
[a,b]) $ (see figure \ref{fig:OD}).
\end{thm}
\begin{proof}
Let us define two maps $\alpha , \beta : I \rightarrow I$ by the
formulas
$$\alpha(x) = F(a,x)$$
$$\beta(x) = F(x,b).$$
Note that $\alpha$ and $\beta$ form something that looks like (but
is not exactly) a $\mathcal{P}$-configuration in $I$. We consider
the dynamical system $(I,\alpha, \beta)$. By the definitions of
$\alpha, \beta$ and by the conditions on $F$ we have that
$$\alpha(b) = \beta(a)$$
and that
$$\alpha(x_0) = a \, \, \textrm{ and } \, \, \beta(y_0) = b $$
and thus
$$ \alpha(I) \cup \beta(I) = I . $$
In addition
$$\forall x \neq y . |\beta(x)-\beta(y)|,|\alpha(x)-\alpha(y)|<|x-y|$$
so all the conditions of proposition \ref{prop:contracting} are
fulfilled and we conclude that the orbit-set of any point in $I$
is dense in $I$.

Now let $f_1$ and $f_2$ be continuous and satisfy
(\ref{eq:boundcond}) and (\ref{eq:feongamma}). We shall show that
for any $z$ in the orbit-set of $a$
$$ f_1(z) = f_2(z) .$$
For $a$ we already have by (\ref{eq:boundcond}) that
$$f_1(a) = A = f_2(a) .$$
If z is a point for which we know that $f_1(z) = f_2(z)$ then
\begin{displaymath}
f_1(\alpha(z)) = f_1(F(a,z)) = H[f_1(a),f_1(z),a,z]
\end{displaymath}
by (\ref{eq:feongamma}). But by our assumption on $z$ we can
replace $H[f_1(a),f_1(z),a,z]$ by $H[f_2(a),f_2(z),a,z]$ and
obtain
\begin{displaymath}
f_1(\alpha(z)) = H[f_2(a),f_2(z),a,z] = f_2(\alpha(z))
\end{displaymath}
where the last equality follows again from (\ref{eq:feongamma}).
So we have
$$f_1(\alpha(z)) = f_2(\alpha(z)) \,\, .$$
Arguing in just the same manner we arrive at the relation
$$f_1(\beta(z)) = f_2(\beta(z)) \,\,.$$
So all the points in the orbit-set of $a$ inherit from $a$ the
property of being given the same values by $f_1, f_2$, and so
indeed for any $z \in OS (a)$ we have $f_1(z) = f_2(z)$. The
continuity of $f_1 , f_2$ and the density of $OS (a)$ imply $f_1 =
f_2$ on $I$.
\end{proof}
As a corollary of the above theorem we have the overdeterminedness
of Jensen's functional equation.
\begin{cor}
Let $\alpha$ and $\beta$ be two positive numbers satisfying
$\alpha+\beta = 1$, and let $I=[a,b]$ be some closed interval.
Then all continuous solutions $f$ of the functional equation
\begin{displaymath}
f(\alpha x + \beta y) = \alpha f(x) + \beta f(y) \,\, , \,\, (x,y)
\in I^2
\end{displaymath}
are of the form
\begin{displaymath}
f(z) = \lambda z + \mu
\end{displaymath}
for some constants $\lambda, \mu \in \mathbb{R}$. Moreover, these
solutions are already determined by the functional equation
\begin{displaymath}
f(\alpha x + \beta y) = \alpha f(x) + \beta f(y) \,\, , \,\, (x,y)
\in \Gamma
\end{displaymath}
where $ \Gamma = ([a,b] \times \{b\} ) \cup (\{a\} \times [a,b])
$.
\end{cor}
\begin{rem}\emph{
As another example of a functional equation that satisfies the
conditions of the theorem, one may take (on an appropriate
interval) the equation of the geometric mean
\begin{displaymath}
f\left(\sqrt{xy}\right) = \frac{1}{2}f(x) + \frac{1}{2}f(y) \quad
.
\end{displaymath}}
\end{rem}

\begin{rem}\emph{Note that the above proof suggests an
algorithm that can compute numerically a solution (when such
exists) to a given functional equation on an interval with
boundary data.}
\end{rem}

\begin{rem}\emph{ Note that it follows from the above theorem that
usually (\ref{eq:feaczel}) will not have a solution, even if
(\ref{eq:feongamma}) has a solution.}
\end{rem}

\chapter{Boundary value problems for hyperbolic PDE's}\label{ch:BVP}
We shall now give two applications of the results in chapter 2 to
PDE's. In the first section we will prove a theorem stating a
necessary and sufficient condition for the well posed-ness of a
third order, strictly hyperbolic partial differential boundary
value problem in the plane. This condition is stated in terms of
the dynamical behavior of some dynamical system on the boundary of
the problem. In the second section we shall translate this
condition into explicit, sufficient conditions for solvability in
terms of the geometric structure of the boundary of the problem.

\section{Formulation of the problem and main result}
\begin{figure}[t]
\begin{minipage}[t]{0.33\textwidth}
\includegraphics[scale=0.5]{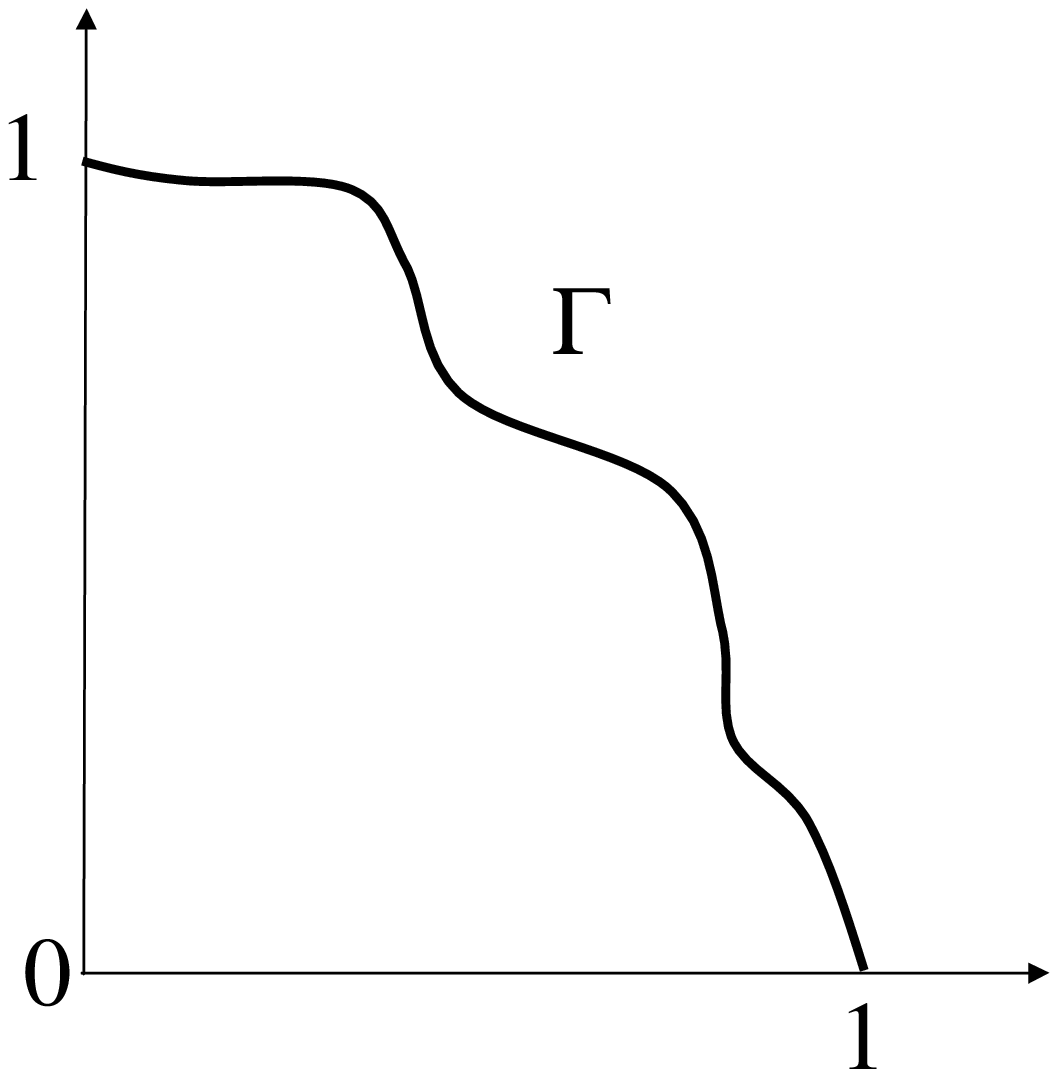}
 \caption{A domain $D$ of the type we consider.}
\label{fig:gooddomain}
\end{minipage}%
\hfill
\begin{minipage}[t]{0.33\textwidth}
\includegraphics[scale=0.5]{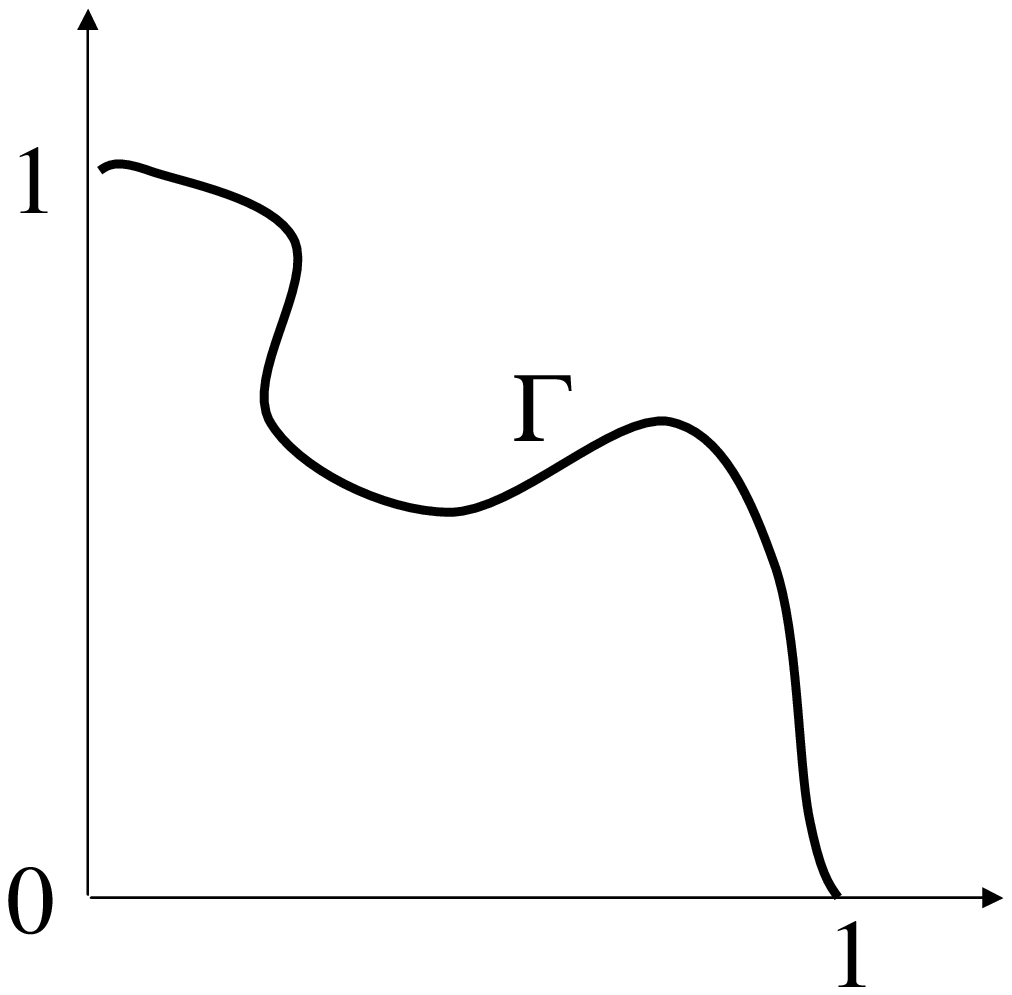}
\caption{A domain $D$ of the type we don't consider.}
\label{fig:baddomain}
\end{minipage}%
\end{figure}

Let $O$ denote the origin in $\mathbb{R}^2$, and let $A_1 = (1,0)$
and $A_2 = (0,1)$. Let $\Gamma$ be a $C^2$ curve that intersects
the axes exactly at the points $A_1$ and $A_2$. We assume that
$\Gamma = \{(\alpha_1(z), \alpha_2(z)) \, | \, z \in [-1,1] \}$
where $\alpha_1, \alpha_2 \in C^2([-1,1])$ satisfy
\begin{equation}\label{eq:monotonealpha}
{\alpha_1}' \geq 0 \quad {\rm and} \quad {\alpha_2}' \leq 0 \,\, .
\end{equation}
We will be dealing with the following problem:
\begin{eqnarray}\label{eq:SPCBP1}
(m\partial_x + n\partial_y)\partial_x \partial_y u &=& 0 \quad {\rm in} \quad D \\
u & = & g \quad {\rm on} \quad \partial D \label{eq:SPCBP2}
\end{eqnarray}
where $m,n>0$ and the domain $D$ is the curvilinear triangle $O
A_1 A_2$ (see figure \ref{fig:gooddomain}) . As for $g$, it is
assumed to be an arbitrary $C^2(\partial D)$ function.
This problem is a special case of the ``second partly
characteristic boundary value problem"\footnote{\cite{Pa05}}.
Under the above assumptions, Paneah proved in \cite{Pa04DAIG}
sufficient conditions for the well-posedness of the problem, which
are also necessary, under some additional assumptions. We shall
exploit the methods introduced in that paper to arrive at a
necessary and sufficient condition for the well-posedness of the
problem under the above assumptions only.

\begin{figure}[t]
\centering
\includegraphics[scale=0.5]{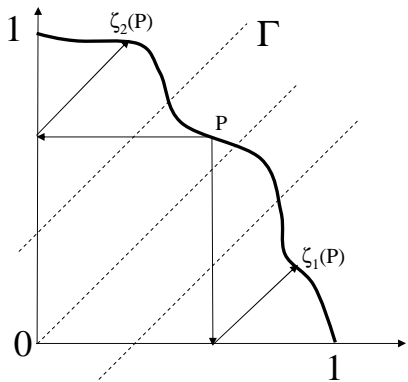}
\caption{The guided dynamical system defined on $\Gamma$.}
\label{fig:gds}
\end{figure}

We construct a guided dynamical system on $\Gamma$. For any point
$p$ in $\overline{D}$, define $\pi_1 p$ to be the projection of
$p$ onto the $x$-axis and $\pi_2 p$ to be the projection of $p$
onto the $y$-axis. Through $p$ there is a line $\ell = \{p +
(mt,nt) : t \in \mathbb{R} \}$. Let $\pi_3 p$ to be the unique
point of intersection of the line $\ell$ passing through $p$ and
of $\Gamma$. Note that $\pi_1$, $\pi_2$ and $\pi_3$ project along
characteristic lines of the operators $\partial y$, $\partial x$
and $m\partial x + n \partial y$, respectively.  We now define two
maps in $\Gamma$
\begin{displaymath}
\zeta_1 = \pi_3 \circ \pi_1 \quad {\rm and} \quad \zeta_2 = \pi_3
\circ \pi_2 \quad .
\end{displaymath}
We introduce the guiding sets
\begin{displaymath}
\Omega_1 = \{q \in \Gamma \, | \, \, (0,1)_q \in T_q (\Gamma)\}
\end{displaymath}
and
\begin{displaymath}
\Omega_2 = \{q \in \Gamma \, | \, \, (1,0)_q \in T_q (\Gamma)\}
\quad .
\end{displaymath}
In words: the set $\Omega_1$ is precisely the subset of $\Gamma$
consisting of points where the tangent line is parallel to the
$y$-axis, and a similar statement holds for $\Omega_2$ (see figure
\ref{fig:gds}). It turns out that the dynamical properties of
$(\Gamma, \zeta, \Omega)$ determine precisely the solvability of
the homogeneous problem (\ref{eq:SPCBP1})-(\ref{eq:SPCBP2}).
Before anything else we must explain what we mean by ``a solution
to (\ref{eq:SPCBP1})-(\ref{eq:SPCBP2})".
\begin{defin}
A function $u \in C(\overline{D})$ is called a \emph{generalized
solution} to the problem (\ref{eq:SPCBP1})-(\ref{eq:SPCBP2}) if
\begin{displaymath}
u \Big| _{\partial D} = g
\end{displaymath}
and if for any $\varphi \in C_0^\infty (D)$
\begin{displaymath}
-\int_D u \, (m\partial_x + n\partial_y)\partial_x \partial_y \,
\varphi \, dxdy = 0 \quad .
\end{displaymath}
\end{defin}

It is convenient to reduce the boundary problem we are interested
in to the problem studied in subsection \ref{subsec:initvalprob}.
Define a map $\omega_{\Gamma}: \Gamma \rightarrow I := [-m,n]$ by
\begin{displaymath}
\omega_{\Gamma}(x,y) = nx - my \,\,.
\end{displaymath}
Define two maps $\delta_1,\delta_2$ in $I$ by
\begin{displaymath}
\delta_i = \omega_{\Gamma} \circ \zeta_i \circ
{\omega_{\Gamma}}^{-1} \,\, .
\end{displaymath}
If we denote $\Lambda_i = \{t \in I \, | \, \delta_i' (t) = 0 \}$,
then lemma 5 in \cite{Pa04DAIG} tells us that
$(\Gamma,\zeta,\Omega)$ and $(I,\delta,\Lambda)$ are isomorphic as
guided dynamical systems. By lemma 4 in that paper
$(I,\delta,\Lambda)$ is a $\mathcal{P}$-configuration. In that
same paper it is shown that finding a generalized solution to
problem (\ref{eq:SPCBP1})-(\ref{eq:SPCBP2}) is \emph{equivalent}
to the problem of finding some function $f \in C^2(I)$ satisfying
the following conditions:
\begin{eqnarray*}
f(t) - f(\delta_1 (t)) - f(\delta_2 (t)) & = & h(t) \quad , \quad
t \in I
\\ \label{eq:GPconfeq2} f'(0) & = & 0
\end{eqnarray*}
where $h$ is an arbitrary $C^2(I)$ function satisfying the
boundary conditions $h(-m)=h(n)$. Combining this reduction and
theorems \ref{thm:initvalprob} and \ref{thm:isomorphic} we
immediately obtain :
\begin{thm}\label{thm:mainthm}
Let $D$ and $\Gamma$ be the domain and the curve described above.
For any $g \in C^2(\Gamma)$ there exists a unique generalized
solution $u \in C^2(\overline{D})$ of the problem
(\ref{eq:SPCBP1})-(\ref{eq:SPCBP2}) if and only if $(\Gamma,
\zeta, \Omega)$ is $\Omega$-minimal.
\end{thm}

\begin{rem}\emph{
It should be noted that in 1941 F. John obtained
results\footnote{\cite{John}} tying the solvability of the
Dirichlet problem for the wave equation with the dynamical
behavior of some system on the boundary of a domain generated by
the characteristics of the wave operator $\partial_x \partial_y$.}
\end{rem}

\section{Explicit conditions for solvability} Theorem
\ref{thm:mainthm} gives the intimate connection between the
dynamical system generated on $\Gamma$ by the characteristic lines
of the operators $\partial x $, $\partial y$ and $m\partial_x +
n\partial_y$ and the unique solvability of the second partly
characteristic boundary value problem. But the condition in the
theorem might seem rather vague. It would be very interesting to
find a geometrical condition on $\Gamma$ that is necessary and
sufficient for $(\Gamma, \zeta, \Omega)$ to be $\Omega$-minimal,
but, unfortunately, we have not been able to find such a
condition. In this section we give explicit conditions that are
sufficient for unique solvability. These conditions will allow us
to ``solve" problem (\ref{eq:SPCBP1})-(\ref{eq:SPCBP2}) in domains
that couldn't be dealt with within the framework of the theory
developed until now. In this section we stick with the notation of
the previous section.

Before stating our results, let us review the results already
known.
\begin{defins}Let $\mathcal{O} = (p_1, p_2, \ldots p_N)$ be an
orbit
\begin{itemize}
  \item If all the points in $\mathcal{O}$ belong to $\Omega$ then $\mathcal{O}$ is called an
  \emph{$\Omega$-guided} orbit.
  \item If $p_1 = p_N$ then $\mathcal{O}$ is called a cycle.
\end{itemize}
\end{defins}
We denote by $\mathcal{N}^\Omega_\zeta$ the set of all
$\Omega$-proper, $\Omega$-guided cycles in $\Gamma$. Also, we let
$\Omega_j '$ be the set of limit points of $\Omega_j$, for
$j=1,2$. In \cite{Pa04DAIG} it is proved that if the following two
conditions hold:
\begin{description}
  \item[1]$\Gamma$ is transversal to the $x$ and $y$ axes at $A_1$ and
$A_2$,
  \item[2]All possible pairs of points $p_1 \in \Omega_1 '$
and $p_2 \in \Omega_2 '$ are situated on $\Gamma$ in the order
$A_2,p_1,p_2,A_1$;
\end{description}
then problem (\ref{eq:SPCBP1})-(\ref{eq:SPCBP2}) is uniquely
solvable if and only if $\mathcal{N}^\Omega_\zeta = \emptyset$. We
do not know wether $\mathcal{N}^\Omega_\zeta = \emptyset$ is a
sufficient condition for unique solvability when conditions {\bf
1} and {\bf 2} are not fulfilled\footnote{Necessity remains, since
an $\Omega$-proper, $\Omega$-guided cycle is an
$(\Omega,\zeta)$-invariant closed subset of $\Gamma$. }. Our main
purpose in this section is to prove the solvability of the second
partly characteristic boundary value problem in domains not
satisfying conditions {\bf 1} and {\bf 2}.

\begin{figure}[t]
\centering
\includegraphics[scale=0.5]{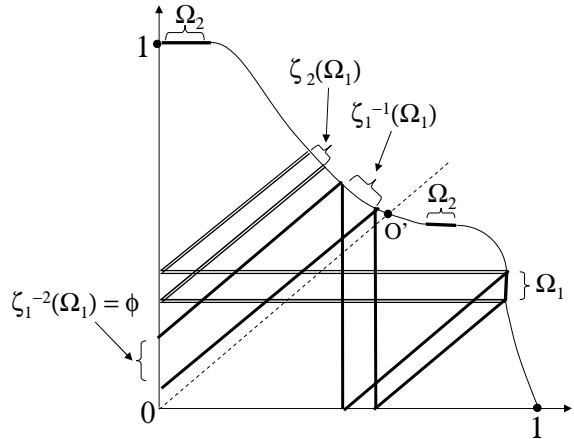}
\caption{Illustration of proposition \ref{prop:exp1}.}
\label{fig:expl1}
\end{figure}

\begin{prop}\label{prop:exp1}
Assume that $A_1 \notin \Omega$, and assume that for any point $p
\in \Omega_1$ there is an $\Omega$-proper orbit
\begin{displaymath}
(p,p_1, \ldots, p_N)
\end{displaymath}
such that \begin{displaymath} p_N \notin \bigcup_{k \geq 0}
{\zeta_1}^{-k} (\Omega_1)
\end{displaymath}
that is, for no $k \geq 0$ the inclusion ${\zeta_1}^k (p_N) \in
\Omega_1$ is possible. Then for any $g \in C^2(\partial D)$ there
exists a unique generalized solution $u \in C^2(\overline{D})$ to
(\ref{eq:SPCBP1})-(\ref{eq:SPCBP2}).
\end{prop}
\begin{proof}
Due to theorem \ref{thm:mainthm}, it is enough to prove that
$(\Gamma,\zeta,\Omega)$ is $\Omega$-minimal.
By proposition \ref{prop:miniffwa} it is enough to prove that
$(\Gamma,\zeta,\Omega)$ has an $\Omega$-weak attractor not in
$\Omega$. We shall show that the point $A_1$ is the desired weak
attractor.

First, note that for any $p \in \Gamma$ the sequence $(p,
\zeta_1(p), \zeta_1^2(p), \ldots )$ converges to $A_1$. Now for
any point $p_0 \in \Gamma$, consider the longest $\Omega$-proper
orbit
\begin{displaymath}
(p_0, \zeta_1 (p_0), \zeta_1^2(p_0), \ldots ) \,\, .
\end{displaymath}
If this orbit is infinite, then it converges to $A_1$ by the above
remark. If this orbit is finite, then there is some $N \geq 0$
such that $(p_0, \zeta_1 (p_0), \zeta_1^2(p_0), \ldots,
\zeta_1^N(p_0))$ is $\Omega$-proper but $\zeta_1^N(t_0) \in
\Omega_1$. But by the assumption of the theorem there is an
$\Omega$-proper orbit $(\zeta_1^N(p_0), p_1, \ldots, p_M)$ such
that $p_M \notin \bigcup_{k \geq 0} {\zeta_1}^{-k} (\Omega_1)$.
Thus
\begin{displaymath}
(p_0, \zeta_1 (p_0), \zeta_1^2(p_0), \ldots, \zeta_1^N(p_0), p_1,
\ldots, p_M, \zeta_1 (p_M), \zeta_1^2 (p_M), \zeta_1^3 (p_M),
\ldots )
\end{displaymath}
is $\Omega$-proper and converges to $A_1$.
\end{proof}
\begin{expl}\label{expl:1}\emph{
Consider the domain in figure \ref{fig:expl1}. Although $\Gamma$
does not satisfy neither of conditions {\bf 1} or {\bf 2} above,
yet, by the above proposition, the boundary problem}
\begin{eqnarray*}
(\partial_x + \partial_y)\partial_x \partial_y u &=& 0 \quad {\rm in} \quad D \\
u & = & g \quad {\rm on} \quad \partial D
\end{eqnarray*}
\emph{has a unique generalized solution $u \in C^2(\overline{D})$
for any $g \in C^2(\partial D)$.}
\end{expl}
The above example is not contained in the results that appeared up
to now, but could have been obtained using the same techniques.
The next proposition deals with a configuration that truly makes
use of the methods and notions that we have introduced in this
thesis.

Before we state our next proposition, let us make some further
notation and remarks. There is some point $z_0 \in [-1,1]$ such
that $(\alpha_1(z_0),\alpha_2(z_0)) = \zeta_1(A_2) =
\zeta_2(A_1)$. Denote $O' = (\alpha_1(z_0),\alpha_2(z_0))$. By
\emph{the open segment $A_2 O'$} we shall mean the homeomorphic
image of $[-1,z_0)$ in $\Gamma$. The open segment $O' A_1$ is
defined similarly.

Now let $(I,\delta,\Lambda)$ be the guided dynamical system
defined in the discussion before theorem \ref{thm:mainthm}. The
map $\delta_1 \circ \delta_2 : I \rightarrow I$ (and similarly,
$\delta_2 \circ \delta_1$) satisfies the inequality
\begin{equation}\label{eq:deltatag1}
\left(\delta_1 \circ \delta_2 \right)'(t) = \delta_2'(t) \cdot
\delta_1'(\delta_2(t)) \leq 1
\end{equation}
for all $t$ in $I$. For all points $t$ not contained in $\Lambda$,
this map satisfies the stronger inequality
\begin{equation}\label{eq:delta_tag2}
\left(\delta_1 \circ \delta_2 \right)'(t) < 1 \,\,.
\end{equation}
The following lemma shows the importance of conditions
(\ref{eq:deltatag1})-(\ref{eq:delta_tag2}).
%
\begin{lem}
Let $I$ be a closed interval, and let $f:I \rightarrow I$ be a
non-decreasing  $C^1$ function satisfying $f' \leq 1$. Let $t_0
\in I$ be a fixed point\footnote{It is well known that every
continuous function from an interval into itself has a fixed
point.} of $f$. If $f'(t_0) < 1$ then $t_0$ is the \emph{unique}
fixed point of $f$ and, moreover, for any $t \in I$ the sequence
\begin{displaymath}
(t, f(t), f^2(t), f^3(t), \ldots)
\end{displaymath}
converges to $t_0$.
\end{lem}
\begin{proof}
Assume, without loss of generality, that $t_0$ is an inner point
of $I$. Define a function $T: I \rightarrow I$ by $T(t) = t$. We
have that $f' \leq T'$ in $I$. Let $t_1>t_0$. Since
$f'(t_0)<1=T'(t_0)$, there is a neighborhood $U$ of $t_0$ such
that $f'(t)<T'(t)$ for all $t \in U$. Thus
\begin{equation*}
f(t_1) - f(t_0)= \int_{t_0}^{t_1}f'(t)\,dt <
\int_{t_0}^{t_1}T'(t)\,dt = t_1 - t_0\,\,.
\end{equation*}
but $t_0$ is a fixed point of $f$, thus
\begin{equation}\label{eq:f(t_1)<t_1}
f(t_1) < t_1 \,\,.
\end{equation}
As a consequence, no point $t_1$ greater than $t_0$ can be a fixed
point of $f$. On the other hand,
\begin{equation}\label{eq:f(t_1) >= t_0}
f(t_1) \geq f(t_0) = t_0
\end{equation}
because $f$ is non-decreasing. The combination of
(\ref{eq:f(t_1)<t_1}) and (\ref{eq:f(t_1) >= t_0}) implies that
the sequence
\begin{equation*}
(t_1, f(t_1), f^2(t_1), f^3(t_1), \ldots)
\end{equation*}
converges to some $t_2 \geq t_0$. Thus $t_2$ is a fixed point of
$f$, and using (\ref{eq:f(t_1)<t_1}) again we conclude that $t_2 =
t_0$. In a similar we may obtain the same results for $t_1 < t_0$.
This completes the proof of the lemma.
\end{proof}

\begin{figure}[t]
\begin{minipage}[t]{0.33\textwidth}
\includegraphics[scale=0.5]{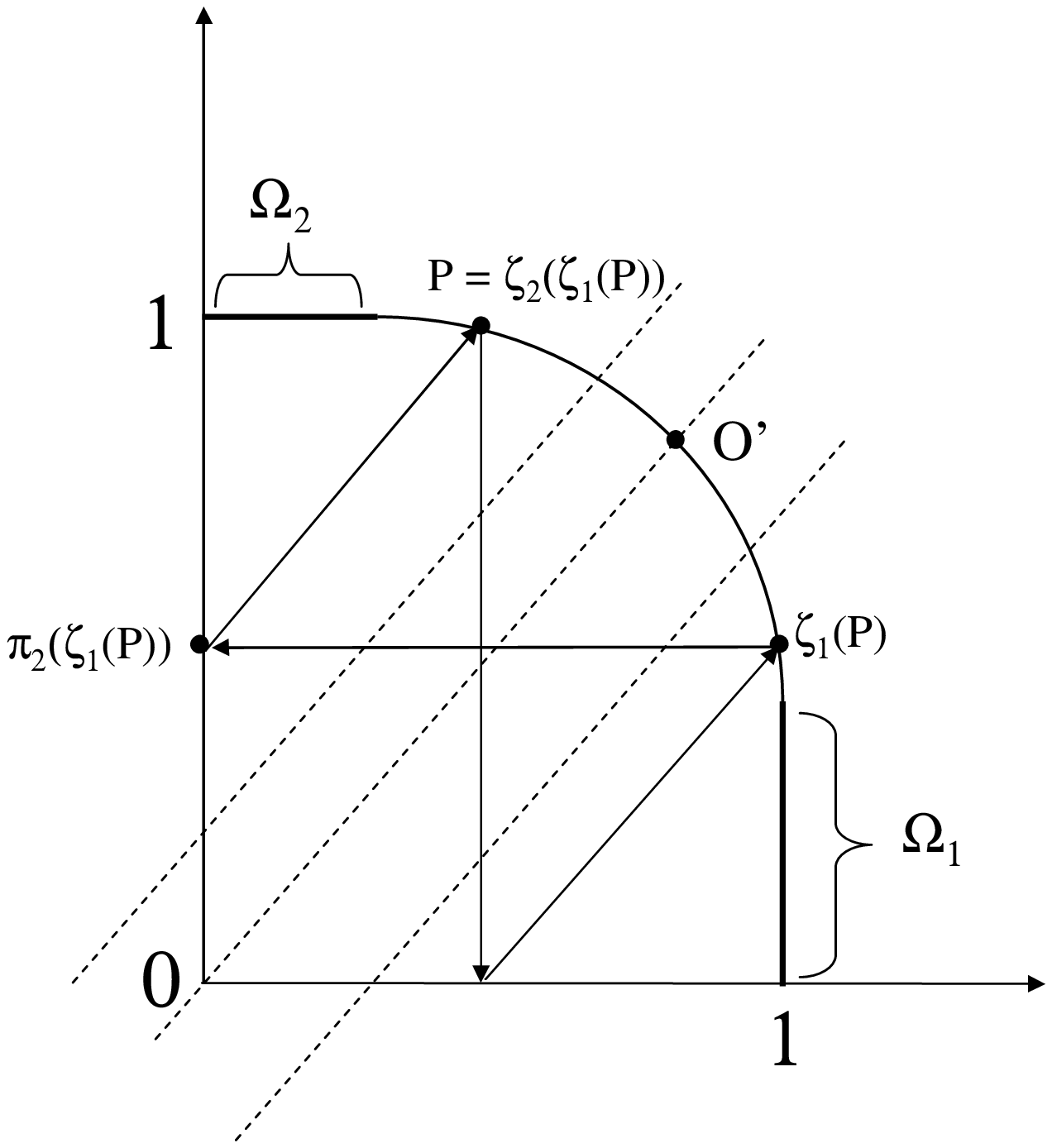}
\caption{Illustration of proposition \ref{prop:exp2}, solvable
case.} \label{fig:expl2}
\end{minipage}%
\hspace{4cm}
\begin{minipage}[t]{0.33\textwidth}
\centering
\includegraphics[scale=0.5]{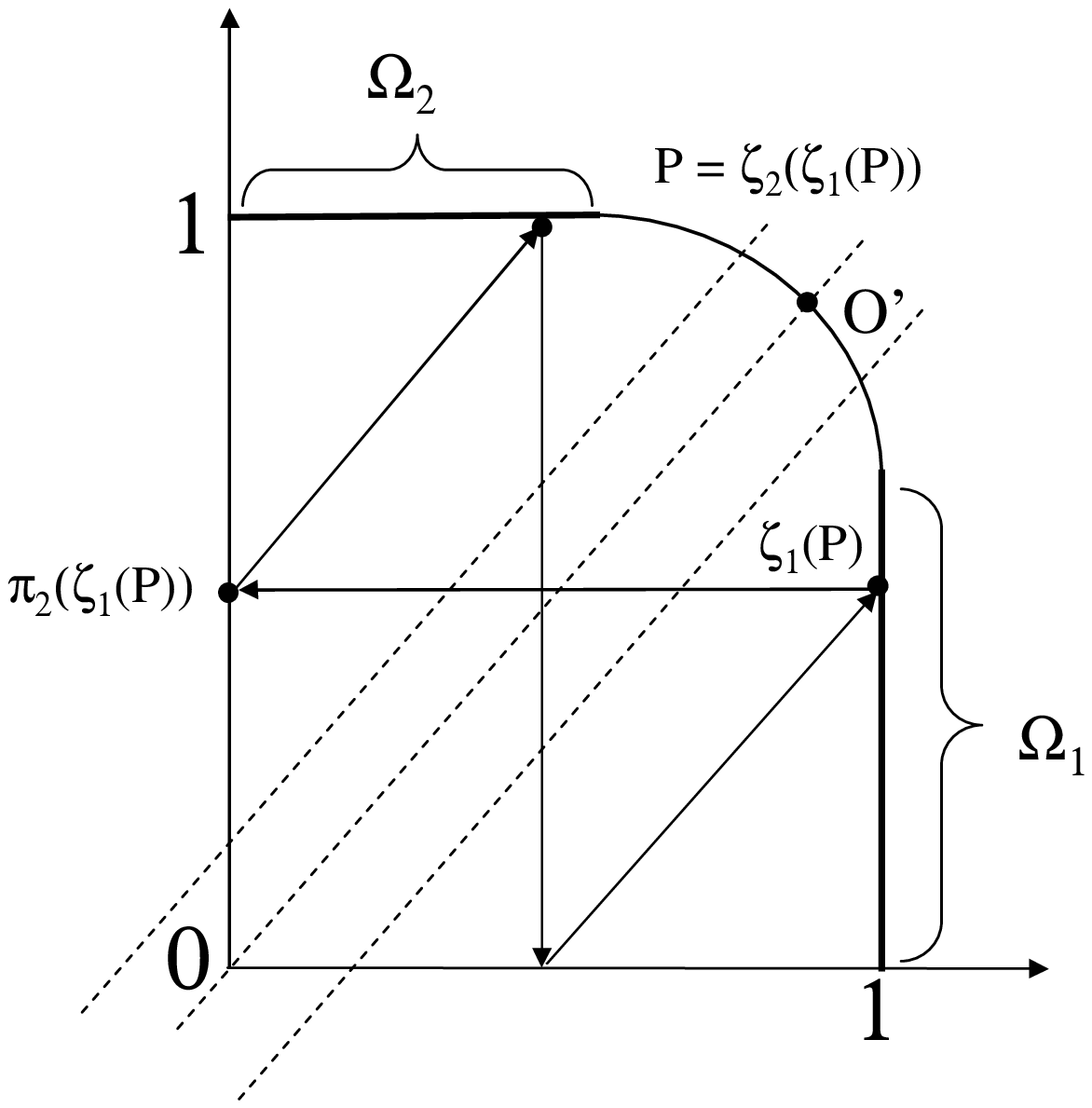}
\caption{Illustration of proposition \ref{prop:exp2}, the
non-solvable case.} \label{fig:expl3}
\end{minipage}%
\end{figure}

Now, assume that the map $\zeta_1 \circ \zeta_2$ has a fixed point
$p_1$ outside of $\Omega$. This fixed point is mapped by the
isomorphism\footnote{See the discussion preceding theorem
\ref{thm:mainthm}.} $w_\Gamma$ to a fixed point $t_1 \notin
\Lambda$ of $\delta_1 \circ \delta_2$. By the discussion before
the lemma,
\begin{displaymath}\label{eq:deltatag2}
\left(\delta_1 \circ \delta_2 \right)'(t_1) < 1 \,\,.
\end{displaymath}
Using the lemma, we conclude that $t_1$ is an attractive fixed
point of $\delta_1 \circ \delta_2$, and this translates to the
fact that for any $p \in \Gamma$, the sequence
\begin{displaymath}
(p, \zeta_1(\zeta_2(p)), \zeta_1(\zeta_2(\zeta_1(\zeta_2(p)))),
\ldots )
\end{displaymath}
converges to $p_1$. The same discussion can also be made for fixed
points of $\zeta_2 \circ \zeta_1$.
\begin{prop}\label{prop:exp2}
Assume that $\Omega_1$ is contained in the open segment $O' A_1$,
and that $\Omega_2$ is contained in the open segment $A_2 O'$.
Then the problem (\ref{eq:SPCBP1})-(\ref{eq:SPCBP2}) has a
(unique) generalized solution $u \in C^2(\overline{D})$ for every
$g \in C^2(\partial D)$ if and only if  either $\zeta_1 \circ
\zeta_2$ or $\zeta_2 \circ \zeta_1$ has a fixed point not in
$\Omega$.
\end{prop}
\begin{proof}
Assume that neither $\zeta_1 \circ \zeta_2$ nor $\zeta_2 \circ
\zeta_1$ have fixed points outside $\Omega$. Denote by $p_1$ a
fixed point of $\zeta_1 \circ \zeta_2$ lying in $\Omega$. Because
$p_1 \in \overline{O' A_1}$, $p_1$ must be in $\Omega_1$. Let $p_2
= \zeta_2(p_1)$. Now,
\begin{displaymath}
\zeta_2(\zeta_1(\zeta_2(p_1))) = \zeta_2(p_1) \,\,,
\end{displaymath}
so $p_2$ is a fixed point of $\zeta_2 \circ \zeta_1$. By
assumption, $p_2 \in \Omega$, but because $p_2 \in \overline{A_2
O'}$, $p_2$ must be in $\Omega_2$. Therefore, set $\{p_1, p_2\}$
is a closed, $(\Omega,\zeta)$-invariant set in $\Gamma$
\footnote{In other words, $\{p_1, p_2\}$ is an $\Omega$-proper,
$\Omega$-guided cycle.}. As $(\Gamma, \zeta, \Omega)$ cannot be
$\Omega$-minimal, theorem \ref{thm:mainthm} tells us that there
are $g \in C^2(\partial D)$ for which there is no solution $u$ to
(\ref{eq:SPCBP1})-(\ref{eq:SPCBP2}).

Now assume, without loss of generality, that a fixed point $p_1$
of $\zeta_1 \circ \zeta_2$ is not in $\Omega$. It is enough to
show that $p_1$ is an $\Omega$-weak attractor. Let $p$ be a point
in the closed segment $O' A_1$. As $\Omega_2$ is contained in the
open segment $A_2 O'$, the orbit $(p, \zeta_2(p))$ is
$\Omega$-proper. Now, $\zeta_2(p)$ is in the closed segment $A_2
O'$, so $(p, \zeta_2(p), \zeta_1(\zeta_2(p)))$ is also an
$\Omega$-proper orbit. Continuing in this fashion, we get an
$\Omega$-proper orbit
\begin{displaymath}
(p, \zeta_2(p), \zeta_1(\zeta_2(p)), \zeta_2(\zeta_1(\zeta_2(p))),
\ldots )
\end{displaymath}
But this orbit contains the sub-sequence
\begin{displaymath}
(p, \zeta_1(\zeta_2(p)), \zeta_1(\zeta_2(\zeta_1(\zeta_2(p)))),
\ldots )
\end{displaymath}
which, as we have mentioned above, converges to $p_1$. So for any
$p$ in the closed segment $O' A_1$, we have $p_1 \in
\overline{\Omega{\rm -}OS(p)}$. Now if $p$ is in the closed
segment $A_2 O'$, then $(p, \zeta_1(p))$ is $\Omega$-proper and
$\zeta_1(p)$ is now in the open segment $O' A_1$. So $p_1 \in
\overline{\Omega{\rm -}OS(\zeta_1(p))}$ which clearly implies that
$p_1 \in \overline{\Omega{\rm -}OS(p)}$. We have established the
fact that $p_1$ is an $\Omega$-weak attractor, and the proof is
completed by calling into action proposition \ref{prop:miniffwa}
and theorem \ref{thm:mainthm}.
\end{proof}
\begin{expl}\label{expl:2}
\emph{The domain in figure \ref{fig:expl2} does not satisfy
condition {\bf 1} nor condition {\bf 2}, but the map $\zeta_1
\circ \zeta_2$ has a fixed point $p \notin \Omega$ so by the above
proposition the boundary value problem}
\begin{eqnarray*}
(\partial_x + \partial_y)\partial_x \partial_y u &=& 0 \quad {\rm in} \quad D \\
u & = & g \quad {\rm on} \quad \partial D
\end{eqnarray*}
\emph{has a unique generalized solution $u \in C^2(\overline{D})$
for any $g \in C^2(\partial D)$.}
\end{expl}

\begin{expl}\label{expl:3}
\emph{Consider figure \ref{fig:expl3}. The fixed points of
$\zeta_1 \circ \zeta_2$ and $\zeta_2 \circ \zeta_1$ are in
$\Omega$, so by the above proposition there are functions $g \in
C^2(\partial D)$ for which the boundary value problem}
\begin{eqnarray*}
(\partial_x + \partial_y)\partial_x \partial_y u &=& 0 \quad {\rm in} \quad D \\
u & = & g \quad {\rm on} \quad \partial D
\end{eqnarray*}
\emph{has no solution $u \in C^2(\overline{D})$.}
\end{expl}

We are able to state the last proposition using the compact
condition that Paneah used. Recall that $\mathcal{N}^\Omega_\zeta$
is the set of all $\Omega$-proper, $\Omega$-guided cycles in
$(\Gamma, \zeta, \Omega)$. Then we have
\begin{cor}
Under the assumptions of proposition \ref{prop:exp2}, the problem
(\ref{eq:SPCBP1})-(\ref{eq:SPCBP2}) has a (unique) generalized
solution $u \in C^2(\overline{D})$ for every $g \in C^2(\partial
D)$ if and only if $\mathcal{N}^\Omega_\zeta = \emptyset$.
\end{cor}
\begin{proof}
We have already mentioned that $\mathcal{N}^\Omega_\zeta =
\emptyset$ is a necessary condition for unique solvability. Now if
$\mathcal{N}^\Omega_\zeta = \emptyset$, the above proof shows that
for either $\zeta_1 \circ \zeta_2$ or $\zeta_2 \circ \zeta_1$
there has to be a fixed point $p \notin \Omega$, thus the above
proposition implies the assertion of the corollary.
\end{proof}

Following the idea of the last proposition, we may obtain
sufficient conditions for solvability in terms of the fixed points
of $\zeta_{i_1}\circ\zeta_{i_2}\circ \cdots \circ \zeta_{i_N}$,
for an arbitrary multi-index $(i_1, i_2, \ldots, i_N)$. But we
shall not write down such theorems. It is the author's belief that
the best kind of progress will be made by finding a simple and
analytical (or geometrical) necessary and sufficient condition for
the guided dynamical system $(\Gamma, \zeta, \Omega)$ to be
minimal. We conclude this thesis with a conjecture in this
direction. We state this conjecture in terms of the second partly
characteristic boundary value problem, although it may also be
viewed as a conjecture regarding the minimality of a
$\mathcal{P}$-configuration.
\begin{conj}
The boundary value problem (\ref{eq:SPCBP1})-(\ref{eq:SPCBP2}) is
uniquely solvable if and only if $\mathcal{N}^\Omega_\zeta =
\emptyset$.
\end{conj}

\chapter{Late introduction}\label{ch:intro}
The main theme of this thesis is the use of \emph{guided dynamical
systems} in problems in functional equations and in partial
differential equations. Most of the problems we deal with were
first studied by Paneah in the papers cited in the text, and
originated from the problem dealt with in \cite{Pa97}. In this
introduction we give a brief overview of the results in this
thesis and survey related known results.

\section{Chapter \ref{ch:dys}}
\emph{Guided dynamical systems} are a generalization of dynamical
systems with several generators. A guided dynamical system is
simply a dynamical system in which each of the generating maps
acts only on a subset of the space. The first guided dynamical
systems appeared in \cite{Pa03NC} and \cite{Pa03SFE}. In these
papers the space was an interval or a curve, and on it acted two
generating maps.

Chapter \ref{ch:dys} is a first step in the development of a
general theory of guided dynamical systems. We develop only the
parts of the theory that are used in other parts of this
thesis\footnote{Theorem \ref{thm:existminimal} is exception to
this rule. It was proved because of its beauty, not its
usefulness.}. It is the author's belief that there is much work
left to be done in the general theory of guided dynamical systems,
which appears both potentially applicable to other parts of
Mathematics and interesting in itself.

In section \ref{sec:definitions} we set the notation for
(non-guided) discrete dynamical systems, define the basic terms
and obtain the first (original) result in this work - proposition
\ref{prop:contracting}. The exposition is influenced by two main
approaches: the first that of Paneah, which is non-standard but
convenient for our uses, and the second is the approach of B.
Hasselblat and A. Katok \cite{HaKatok}. We found it necessary to
introduce the term \emph{weak attractor} since, on the one hand,
the notion this term represents plays a key role in the theory we
develop and, on the other hand, such a notion has not been given a
name in the literature.

Sections \ref{sec:guideddys} and \ref{sec:iso} deal with
generalizing well-known definitions and results from the theory of
discrete dynamical systems to guided dynamical systems. A
particular case of lemma \ref{lem:propiffprop} was proved in lemma
5 of \cite{Pa04DAIG}. In that paper there was an isomorphism of
guided dynamical systems between to specific guided dynamical
systems - one on an interval and one on a curve.

\section{Chapter \ref{ch:FuncEqs}}
This chapter is devoted to uniqueness and solvability of
functional equations that have the form
\begin{equation}\label{eq:CTFE}
f(x) - \sum_{i=1}^N a_i (x) f(\delta_i (x)) = h(x) \quad, \,\, x
\in X \,\,.
\end{equation}
Here the functions $a_i$, $\delta_i$ and $h$ are given, and $f$ is
an unknown function on $X$. In \cite{Pa03SFE} Paneah studies this
equation where $X$ was an interval and the $\delta_i$'s were
non-decreasing maps satisfying some conditions.

In subsection \ref{subsec:maxp} we generalize the first parts of
theorems 1 and 2 in \cite{Pa03SFE} in a few directions. Theorem 2
from \cite{Pa03SFE} (the \emph{maximum principle}), originally
stated for continuous functions, is generalized to semi-continuous
functions (lemma \ref{lem:maxp}). Theorem 2 (uniqueness of
solutions), originally stated for scalar valued functions, is
generalized to vector valued functions (\ref{thm:maxp}). In both
cases our results hold for (at least) a general compact metric
space $X$. Our sufficient condition for uniqueness up to an
additive constant is given in terms of the existence of a
$\Lambda$-weak attractor.

Subsection \ref{subsec:CauchyRn} deals with uniqueness of
continuously differentiable solutions of (\ref{eq:CTFE}) with
$N=2$, $a_1 \equiv a_2 \equiv 1$ and $X$ a subset of
$\mathbb{R}^n$. Theorems \ref{thm:cauchyveceq} and
\ref{thm:affineCTFE} are the main results of this subsection. Such
results were obtained by Paneah in \cite{Pa04OD} for classes of
functions defined on an interval and differentiable either at the
origin or on the entire interval - depending on the behavior of
the $\delta_i$'s. As above, our sufficient condition for
uniqueness up to an multiplicative constant is given in terms of
the existence of a $\Lambda$-weak attractor.

The main results in section \ref{sec:UniSovCTFE} is theorem
\ref{thm:uniquesolvability}. This theorem is a generalization of
the second part of theorem 3 in \cite{Pa03NC} (regarding unique
solvability of equation \ref{eq:CTFE}) originally stated for maps
on an interval with an attractor in the boundary of the interval,
to general guided dynamical systems with some weak attractor. We
give a proof that is based on the proof in \cite{Pa03NC}, adding a
few details. Our sufficient condition for unique-solvability is
given in terms of the existence of a $\Lambda$-weak attractor.

In section \ref{subsec:initvalprob} we treat the problem of
existence and uniqueness of $C^2$ solutions $f$ to the problem
\begin{eqnarray}\label{eq:GP1}
f(t) - \sum_{i=1}^N f(\delta_i (t)) & = & h(t) \quad , \quad t \in
I
\\ \label{eq:GP2} f'(c) & = & \mu \quad ,
\end{eqnarray}
where $I = [a,b]$ is an interval, $h \in C^2(I)$ is given and the
maps $\delta_i$, $i=1, \ldots, N$ satisfy what we call a
\emph{generalized}
$\mathcal{P}$-\emph{configuration}\footnote{$\mathcal{P}$-configurations
were introduced by Paneah in the papers cited here.}. To be
precise, we give a necessary and sufficient condition for the
existence and uniqueness of solution to
(\ref{eq:GP1})-(\ref{eq:GP2}). This necessary and sufficient
condition is given in terms of a dynamical property of the guided
dynamical system generated in $I$ by the maps $\delta_i$, $i=1,
\ldots N$. Theorem \ref{thm:initvalprob} states that it is the
$\Lambda$-minimality of this dynamical system that is necessary
and sufficient for the unique solvability of the above problem.
This is an improvement on theorem 9 from \cite{Pa03SFE} for two
reasons. First, the passage from $2$ to $N$ maps is not completely
trivial. Second, we give a necessary and sufficient condition,
whereas until now a necessary and sufficient condition for unique
solvability was known only under some additional conditions.
Another interesting new result in this section is proposition
\ref{prop:miniffwa}, which states that for a generalized
$\mathcal{P}$-configuration the existence of a $\Lambda$-weak
attractor (in some set) is equivalent to $\Lambda$-minimality.

\section{Chapter \ref{ch:OD}}
In chapter \ref{ch:OD} we give two results regarding
\emph{overdeterminedness} of functional equations. Details about
this subject are given in the beginning of that chapter. In
theorem \ref{thm:ODofCE} we prove that if a continuous function
$f:[0,1] \rightarrow \mathbb{R}$ satisfies Cauchy's functional
equation
\begin{equation}\label{eq:CE}
f(x+y) = f(x) + f(y)
\end{equation}
on part of the boundary of the square ${\bf K} = \{(x,y) : |x| +
|y| \leq 1\}$ then $f(z) = \lambda z$. This result is weaker in
some sense than a known result of Lackovich, who showed in
\cite{Lackovich} that if $f: [0,\infty) \rightarrow \mathbb{R}$ is
\emph{measurable} and satisfies Cauchy's equation on the line
$\{(at, bt) | t \in \mathbb{R} \}$ where $\log_a b \notin
\mathbb{Q}$, then $f(z) = \lambda z$. Our contributions are that
we have shown that the continuous solutions of equation
(\ref{eq:CE}) are determined on the boundary of ${\bf K}$, and
also that we give a very simple proof based on dynamical
systems\footnote{Lackovich' proof relied on the Krein-Milman
theorem, whereas our proof uses quite elementary analysis}. One
may also compare our results to theorem 1 in \cite{Pa04OD}, where
it is shown that the $C^1$ solutions of (\ref{eq:CE}) are
determined on even a smaller part of the boundary of ${\bf K}$.

Theorem \ref{thm:uniqueOD} can be viewed either as a uniqueness
theorem or as an overdeterminedness theorem. As a uniqueness
theorem, it is very similar to Acz\'{e}l's theorem
\ref{thm:aczel}, just the conditions are slightly different and
the proof is completely different. As an overdeterminedness
theorem, it is probably the first of its kind.

\section{Chapter \ref{ch:BVP}}
In \cite{Pa04DAIG} Paneah reduces the so-called ``second partly
characteristic" third order strictly hyperbolic boundary value
problem to a Cauchy type functional equation. This reduction,
together with theorems \ref{thm:isomorphic} and
\ref{thm:initvalprob}, immediately imply theorem
\ref{thm:mainthm}. This theorem gives the precise connection
between the dynamics in the boundary of the problem and the
solvability of that problem.

Propositions \ref{prop:exp1} and \ref{prop:exp2} give sufficient
conditions for solvability in domains for which there was no
previous result.


%
%

\bibliographystyle{abbrv}
\clearpage \addcontentsline{toc}{chapter}{Bibliography}
\bibliography{thesis_orr}

\end{document}